\documentclass[12pt]{article}
\usepackage{epsfig}
\usepackage{algorithm,algorithmic}
\usepackage{amsmath}
\usepackage{amssymb}
\usepackage{url}
\DeclareMathOperator{\sign}{sign}
\begin{document}
\title{\vspace{-17mm}
On Computing Elastic Shape Distances between
Curves in $d-$dimensional Space}
%\author{Usual suspects (Javier$^1$, Jim$^{1,2}$, Gunay$^1$, Charles$^1$)\\
\author{\tt\small Javier Bernal$^1$, Jim Lawrence$^{1,2}$, Gunay Dogan$^1$, Charles Hagwood$^1$\\
$^1${\small \sl National Institute of Standards and Technology,} \\
{\small \sl Gaithersburg, MD 20899, USA} \\
$^2${\small \sl George Mason University,} \\
{\small \sl 4400 University Dr, Fairfax, VA 22030, USA} \\
{\tt\small $\{$javier.bernal,james.lawrence,gunay.dogan,charles.hagwood$\}$} \\
{\tt\small @nist.gov \ \ \  lawrence@gmu.edu}}
\date{\ }
\maketitle
\vspace{-13mm}
\begin{abstract}
The computation of the elastic registration of two simple curves in higher dimensions
and therefore of the elastic shape distance between them has been investigated by
Srivastava et al. Assuming the first curve has one or more starting points, and the
second curve has only one, they accomplish the computation, one starting point of
the first curve at a time, by minimizing an $L^2$ type distance between them based
on alternating computations of optimal diffeomorphisms of the unit interval and
optimal rotation matrices that reparametrize and rotate, respectively, one of the
curves. We recreate the work by Srivastava et al., but in contrast to it, again for
curves in any dimension, we present a Dynamic Programming algorithm for computing
optimal diffeomorphisms that is linear, and justify in a purely algebraic manner
the usual algorithm for computing optimal rotation matrices, the Kabsch-Umeyama
algorithm, which is based on the computation of the singular value decomposition of
a~matrix. In addition, we minimize the $L^2$ type distance with a procedure that
alternates computations of optimal diffeomorphisms with successive computations
of optimal rotation matrices for all starting points of the first curve. Carrying
out computations this way is not only more efficient all by itself, but, if both
curves are closed, allows applications of the Fast Fourier Transform for computing
successively in an even more efficient manner, optimal rotation matrices for all
starting points of the first curve.
\\[0.2cm]
\end{abstract}
\section{\large Introduction}
In this paper, following ideas in~\cite{srivastava, srivastava2},
we address the problem of computing the elastic shape distance between two simple
curves, not necessarily closed, in $d-$dimensional space, $d$ a positive integer,
or equivalently the problem of computing the elastic registration of two such curves.
If neither curve is closed so that each curve has a fixed starting point,
this is done by finding a diffeomorphism of the unit interval, and
a $d\times d$ rotation matrix, that reparametrizes and rotates, respectively, one of
the curves, not necessarily the same curve for both operations, so that an $L^2$ type
distance between the curves is minimized. The resulting minimum distance is then the
elastic shape distance between the curves and the result of the reparametrization and
rotation of the curves is their elastic registration. On the other hand, if, say, the
first curve is closed so that any point in it can be treated as a starting point of
the curve, then a finite subset of consecutive points is selected in the curve in the
direction in which the curve is defined, in such a way that the subset is reasonably
dense in the curve (by joining consecutive points in the subset with line segments, the
resulting piecewise linear curve should be a reasonable approximation of the curve). This
finite subset of the first curve is interpreted to be the set of starting points of the
curve. A fixed starting point is then identified on the other curve, the second curve,
perhaps arbitrarily if the curve is closed. In \cite{srivastava,srivastava2}, given a
point in the set above interpreted to be the set of starting points of the first curve,
using the point as the starting point of the (first) curve, an optimal diffeomorphism and
an optimal rotation matrix are found in the same manner as described above for the case
in which neither curve is closed.  Again in \cite{srivastava,srivastava2}, this is done
for each point in the set, and the point for which the $L^2$ type distance is the
smallest is then considered to be the optimal starting point in the first curve, and the
optimal diffeomorphism and optimal rotation matrix associated with it are then treated as
the optimal diffeomorphism and optimal rotation matrix that produce the elastic
registration of the two curves and the elastic shape distance between~them.
\par We note that above we have tacitly assumed that the curves are defined in the proper
directions for the purpose of unambiguously comparing their shapes. We have done this for
simplicity as in reality it may not be the case. Clearly, given a simple
curve in $d-$dimensional space, it has two possible directions in which it can be defined.
In particular, in $2-$dimensional space (the plane), a closed simple curve is defined in
either the clockwise direction or the counterclockwise direction, and if the shapes of two
closed simple curves in the plane are to be compared, it only makes sense that both be
defined in the same direction (clockwise or counterclockwise). Unfortunately, in general,
defining two curves in $d-$dimensional space in the proper directions cannot be done this way,
and the only alternative is first to compute the elastic shape distance and registration with
the curves as given and then reverse the direction of one of the curves and do the computations
again. The smaller of the two computed elastic shape distances then determines the proper
directions of the curves, and therefore their correct elastic registration.
Again for simplicity, in the rest of this paper, given two simple
curves in $d-$dimensional space, we assume they are defined in the proper directions for all
purposes, keeping in mind that if this is not the case, all that is required to fix them,
is that the direction of one of the curves be~reversed.
\par
\begin{figure}
\centering
\begin{tabular}{cc}
\includegraphics[width=0.4\textwidth]{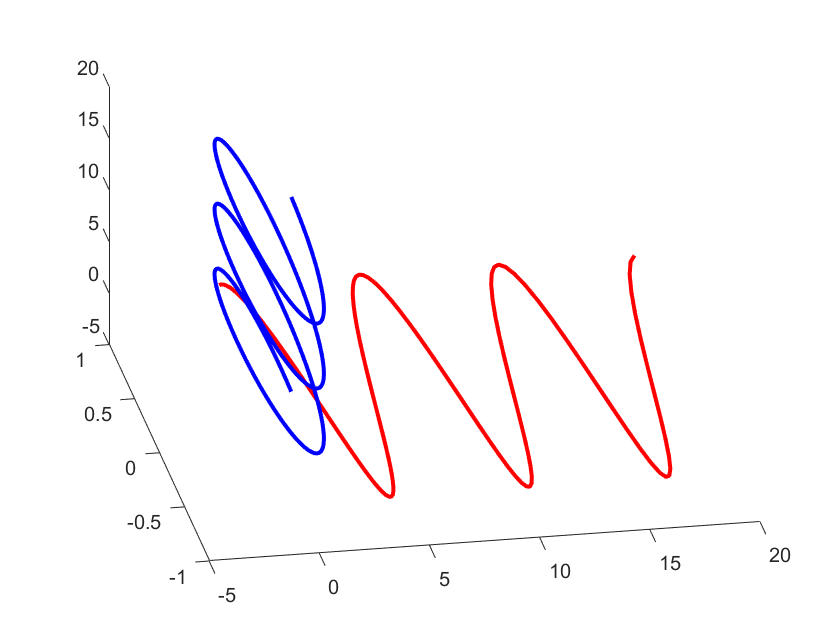}
&
\includegraphics[width=0.4\textwidth]{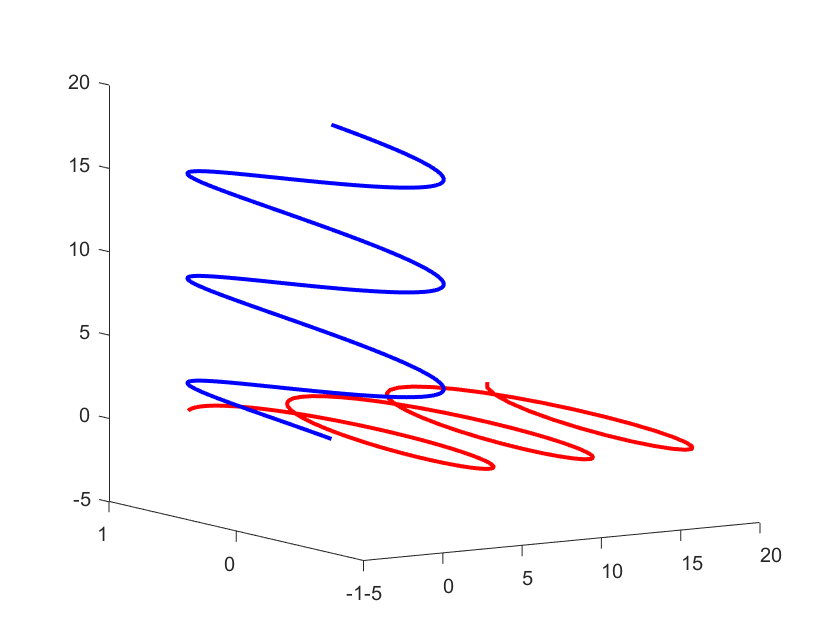}
\end{tabular}
\caption{\label{F:curves0}
Two views of the same two helices, curves in $3-$d space. The positive $z-$axis is the
axis of rotation of one helix, while the positive $x-$axis is the axis of rotation of the other one.
Their shapes are essentially identical thus the elastic shape distance between them should be
essentially~zero.
}

\end{figure}

 Being able to compute the elastic registration of two curves and the elastic
shape distance between them in higher dimensions, in particular in three dimensions,
is useful for studying fiber tracts, geological terrains, protein structures, facial
surfaces, color images, cylindrical helices, etc. See Figure~\ref{F:curves0} that
depicts two such curves in $3-$dimensional space. (Note that in the plots
there, the $y-$axis is not to scale relative to the $x-$axis and the $z-$axis).
\par As mentioned above, in this paper, we address the problem of computing the elastic
registration of two simple curves in $d-$dimensional space, and the elastic shape distance
between them following ideas in~\cite{srivastava,srivastava2}. However, in contrast
to~\cite{srivastava,srivastava2}, we present a method for computing optimal diffeomorphisms
that is linear, and justify in a purely algebraic manner the usual algorithm for computing
optimal rotation matrices. With the convention that if at least one of the curves is closed,
the first curve is closed, so that any point in it can then be treated as a starting point of
the curve, we redefine the $L^2$ type distance to allow for the second curve to be reparametrized
while the first one is rotated, and select, in the appropriate manner, just as it is done in
\cite{srivastava,srivastava2}, a finite subset of points in the (first) curve (one point if
neither curve is closed) which we interpret to be the set of starting points of the curve.
We then minimize the redefined $L^2$ type distance with an iterative procedure that alternates
computations of optimal diffeomorphisms (a constant number of them per iteration for reparametrizing
the second curve) with successive computations of optimal rotation matrices (for rotating the first
curve) for all starting points of the first curve. (Note that in \cite{srivastava,srivastava2}
the alternating computations occur one starting point of the first curve at a time). As noted
in \cite{dogan2}, carrying out computations this way is not only more efficient all by itself,
but, if both curves are closed, allows applications of the Fast Fourier Transform (FFT) as
demonstrated in \cite{dogan} for~$d=2$, for computing successively in an even more
efficient manner, optimal rotation matrices for all starting points of the first~curve.
\par In Section 2 of this paper, we describe a fast linear Dynamic Programming (DP)
algorithm for computing an approximately optimal diffeomorphism for the elastic registration
of two simple curves in $d-$dimensional space, the curves not necessarily closed, each curve
with a fixed starting point, the computation of the registration based only on reparametrizations
(with diffeomorphisms of the unit interval) of one of the curves. In Section~3, we describe and
justify the usual algorithm, the Kabsch-Umeyama algorithm, for computing an approximately optimal
$d\times d$ rotation matrix for the rigid alignment of two simple curves in $d-$dimensional
space, the curves not necessarily closed, each curve with a fixed starting point, the
computation of the alignment based only on rotations (with rotation matrices) of one of the
curves. The algorithm, which is based on the computation of the singular value decomposition
of a matrix, is justified in a purely algebraic manner. In Section~4, given two simple curves
in $d-$dimensional space, one considered to be the first curve, the other the second curve,
and then keeping in mind that the first curve can have one or more starting points while
the second curve has only one, we redefine the $L^2$ type distance between them as hinted
above, and present the iterative procedure mentioned above that minimizes this distance by
alternating computations of approximately optimal diffeomorphisms (a constant number of them
per iteration) with successive computations of approximately optimal rotation matrices for all
starting points of the first curve. In Section~5, assuming both curves are closed, we show
how the FFT can be used to speed up the successive computations of approximately optimal
rotation matrices for all starting points of the first curve.  Finally, in Section~6,
we present results computed with implementations of our methods applied on
%$3-$dimensional
curves in $3-$d space of the helix and spherical ellipsoid~kind.
\section{\large Computation of Diffeomorphism for Registration of Curves in
$d-$dimensional Space}
In this section, we describe a fast Dynamic Programming (DP) algorithm that runs in linear
time to compute an approximately optimal diffeomorphism for the elastic registration of two simple
curves in $d-$dimensional space, $d$ a positive integer, the curves not necessarily closed, each
curve with a fixed starting point. Here the computation of the registration is
based only on reparametrizations (with diffeomorphisms) of one of the curves. Because the
algorithm depends on a couple of input parameters that should be small but sometimes are
too small, it is not guaranteed to succeed every time, but in our experiments we have
observed very convincing results for every problem on which the algorithm was applied,
a few times after adjustments to the parameters.
Here and in what follows, given $T>0$ and an integer~$N>0$, we say a curve is discretized
by a partition $\{t_i\}_{i=1}^{N}$, $t_1=0<t_2<\ldots<t_{N}=T$, of~$[0,T]$, if a
function $\beta: [0,T] \rightarrow \mathbb{R}^d$ has been identified to represent the curve
($\beta$ is continuous and satisfies that its range is exactly the curve), and the curve is
given as the set of $N$ nodes equal to $\{\beta(t_i),i=1,\ldots,N\}$.
On input, the two curves under consideration are given as discrete
sets of nodes in the curves, the result of discretizing the curves by partitions not
necessarily uniform of~$[0,1]$, the numbers of nodes in the curves not necessarily equal.
Given that the numbers of nodes in the
curves are $N$ and $M$, respectively, then the algorithm is indeed linear as it runs in $O(N+M)$
time (see below).
%time as we will show below.
We note that what follows in this section about the algorithm, already appears
in~\cite{bernal} for~$d=2$. We repeat it not only for self-containment but for clarity as
it is what must be said about the algorithm for any $d$ besides~$d=2$.
%We note that the description of the algorithm presented here
%already appears in~\cite{bernal} for~$d=2$.
\par Assume the curves can be represented by functions
$\beta_n: [0,1] \rightarrow \mathbb{R}^d,\ n=1,2$, that are absolutely continuous
(see \cite{bernal2,srivastava}) and of unit length, and given $n$, $n=1,2$, assume
$\beta_n(0)$ is the fixed starting point of $\beta_n$, and if $\beta_n$ is closed, then
$\beta_n(0)=\beta_n(1)$, $\dot{\beta_n}(0)=\dot{\beta_n}(1)$ (here and in what follows,
we may refer to the curve that $\beta_n$ represents simply by~$\beta_n$). With $\|\cdot\|$
as the $d-$dimensional Euclidean norm, we define $q_n: [0,1] \rightarrow \mathbb{R}^d,\ n=1,2$,
by $q_n(t) = \dot{\beta_n}(t) / \Vert\dot{\beta_n}(t)\Vert^{1/2}$ ($d-$dimensional~0 if
$\dot{\beta_n}(t)$ equals $d-$dimensional~0), $q_n$ the shape function or square-root
velocity function (SRVF) of $\beta_n$, and note that because $\beta_n$ is of unit length,
then $\int_0^1 \|q_n(t)\|^2 dt = 1< \infty$, the integral here and in what follows
computed as a Lebesgue integral, so that $q_n$ is square integrable for
$n=1,2$ (see \cite{bernal2,srivastava}). Note as well that for $d=1$,
$q_n(t)$ equals $\sign(\dot{\beta_n}(t))\sqrt{|\dot{\beta_n}(t)|}$, $n=1,2$, the square-root
slope function (SRSF) of~$\beta_n$, and for any $d$, any real number~$C$ and any
square-integrable $q:[0,1]\rightarrow \mathbb{R}^d$ with $\int_0^1 \|q(t)\|^2 dt = 1$, the
function $\beta:[0,1]\rightarrow \mathbb{R}^d$ defined by $\beta(t)=C+\int_0^t q(s)|q(s)|ds$
is absolutely continuous and of unit length with SRVF equal to $q$ almost everywhere
on~$[0,1]$ (see \cite{bernal2,srivastava}).
Finally, we note, given an absolutely continuous function $\beta:[0,1]\rightarrow\mathbb{R}^d$
and $\gamma$ a diffeomorphism of $[0,1]$ onto itself with
$\gamma(0)=0, \gamma(1)=1, \dot{\gamma} > 0$, then $(\beta\circ\gamma)(0)=\beta(0)$, and if
$q$ is the SRVF of $\beta$, then the SRVF of $\beta\circ\gamma$ equals
$(q\circ\gamma)\sqrt{\dot{\gamma}}$ almost everywhere on $[0,1]$ (see \cite{bernal2,srivastava}).
Defining $F(t,\gamma(t),\dot{\gamma}(t)) = \|q_1(t) - \sqrt{\dot{\gamma}(t)} q_2(\gamma(t)) \|^2$,
$\gamma$ as above, we minimize the following energy with respect to~$\gamma$
\begin{equation}\label{E:energy-diff}
E(\gamma) = \int_0^1 F(t,\gamma(t),\dot{\gamma}(t)) dt.
\end{equation}
\par In practice, we need to solve a discretized version of the problem.
Thus, for positive integers $N$, $M$, not necessarily equal, and partitions
of~$[0,1]$, $\{t_i\}_{i=1}^{N}$, $t_1=0<t_2<\ldots<t_{N}=1$, $\{z_j\}_{j=1}^{M}$,
$z_1=0<z_2<\ldots<z_{M}=1$, not necessarily uniform, we assume $\beta_1$ and $\beta_2$
are given as lists of $N$ and $M$ points or nodes in the curves, respectively,
where for $i=1,\ldots,N$, $\beta_1(t_i)$ is the $i^{th}$ point in the list for $\beta_1$,
and for $j=1,\ldots,M$, $\beta_2(z_j)$ is the $j^{th}$ point in the list for $\beta_2$.
We note that $\{\beta_1(t_i), i=1,\ldots,N\}$ is then a set of consecutive points in $\beta_1$
in the direction in which $\beta_1$ is defined, and it should be defined in such a way
that it is reasonably dense in $\beta_1$ (by joining consecutive points in it with
line segments, the resulting piecewise linear curve should be a reasonable approximation
of~$\beta_1$). Similarly for $\{\beta_2(z_j), j=1,\ldots,M\}$ with respect to~$\beta_2$.
We also assume $\dot{\beta_1}(t_i)$, $i=1,\ldots,N$, and $\dot{\beta_2}(z_j)$,
$j=1,\ldots,M$, are approximately computed with centered finite differences from the
lists of points for $\beta_1$ and $\beta_2$, respectively, so  that
$q_1(t_i)$, $i=1,\ldots,N$, and $q_2(z_j)$, $j=1,\ldots,M$, are then approximately
computed by setting $q_1(t_i) = \dot{\beta_1}(t_i) / \Vert\dot{\beta_1}(t_i)\Vert^{1/2}$
($d-$dimensional 0 if $\dot{\beta_1}(t_i)$ equals $d-$dimensional 0), $i=1,\ldots,N$, and
$q_2(z_j) = \dot{\beta_2}(z_j) / \Vert\dot{\beta_2}(z_j)\Vert^{1/2}$
($d-$dimensional 0 if $\dot{\beta_2}(z_j)$ equals $d-$dimensional 0), $j=1,\ldots,M$.
Finally, given $\gamma$, treating \eqref{E:energy-diff} as a Riemann integral, we
discretize \eqref{E:energy-diff} with the trapezoidal rule:
\begin{equation}\label{E:energy-discr}
E(\vec{\gamma}) =\frac{1}{2} \sum_{i=1}^{N-1}
h_i(F(t_{i+1},\gamma_{i+1},\dot{\gamma}_{i+1}) +
F(t_i,\gamma_i,\dot{\gamma}_i)),
\end{equation}
where $h_i = t_{i+1}-t_i$ for $i=1,\ldots,N-1$,
$F(t_i,\gamma_i,\dot{\gamma}_i)$ is understood to be
$\|q_1(t_i) - \sqrt{\dot{\gamma}_i} q_2(\gamma_i) \|^2$ for $i=1,\ldots,N$,
$\vec{\gamma}=(\gamma_i)_{i=1}^N$,
$\gamma_1=0, \gamma_N=1$, $\gamma_i = \gamma(t_i)$,
$\dot{\gamma}_i=(\gamma_{i+1} - \gamma_i)/h_i$ for $i=1,\ldots,N-1$,
$\dot{\gamma}_N=\dot{\gamma}_1$, and $q_2(\gamma_i)$, $i=1,\ldots,N$,
are approximations of $q_2$ at each $\gamma_i$ obtained
from the interpolation of $q_2(z_j)$, $j=1,\ldots,M$, by a cubic~spline. Thus, the
problem of minimizing \eqref{E:energy-diff} with respect to $\gamma$ then becomes,
in practice, the problem of minimizing \eqref{E:energy-discr} with respect to a
discretized~$\gamma$.
\par For this purpose, we consider the $N\times M$ grid on the unit square with grid points
labeled $(i,j)$, $i$, $j$ integers, $1\leq i\leq N$, $1\leq j\leq M$,
each grid point $(i,j)$ coinciding with the planar point~$(t_i,z_j)$.
\par If the mesh of each partition, i.e., $\max (t_{m+1}-t_m), 1\leq m\leq N-1$,
and $\max (z_{m+1}-z_m), 1\leq m\leq M-1$, is sufficiently small, then the set
of diffeomorphisms $\gamma$ of $[0,1]$ onto itself with
$\gamma(0)~=~0, \gamma(1)=1, \dot{\gamma} > 0$, can be approximated by the set of
homeomorphisms of $[0,1]$ onto itself whose graphs are piecewise linear paths from
grid point $(1,1)$ to grid point $(N,M)$ with grid points as vertices, each linear
component of a path having positive slope.  We refer to the latter set as~$\Gamma$.
Then $\gamma$ in~$\Gamma$ is an approximate diffeomorphism of $[0,1]$ onto itself
and as such an energy conceptually faithful to~\eqref{E:energy-discr} can be defined
and computed for it. This is done one linear component of the graph of~$\gamma$
at~a~time.
\par Accordingly, given grid points $(k,l)$, $(i,j)$, $k<i$, $l<j$, that are endpoints
of a linear component of the graph of~$\gamma$, an energy of a trapezoidal nature
over the line segment joining $(k,l)$ and~$(i,j)$ is defined as follows:
%
%\medskip \par
%(\|\hat{q}_1(t'_{m+1})-$
\begin{equation}\label{E:energy6}
%\sqrt{L}q_2(\alpha(t'_{m+1}))\|^2 +
%\|\hat{q}_1(t'_m)-\sqrt{L}q_2(\alpha(t'_m))\|^2).
%(F(t'_{m+1},\alpha(t'_{m+1}),L) +
%F(t'_m,\alpha(t'_m),L).
E_{(k,l)}^{(i,j)}\equiv
\frac{1}{2}\sum_{m=k}^{i-1}(t_{m+1}-t_m)(F_{m+1}+F_m),
\end{equation}
\begin{center}
$F_m\equiv F(t_m,\alpha(t_m),L)$, $m = k,\ldots,i$,
\end{center}
where $F(t_m,\alpha(t_m),L)$ is understood to be 
$\|q_1(t_m) - \sqrt{L} q_2(\alpha(t_m)) \|^2$ for $m = k,\ldots,i$,
where $\alpha$ is the linear function from $[t_k,t_i]$ onto $[z_l,z_j]$ whose
graph is the line segment, $\alpha(t_k)=z_l$, $\alpha(t_i)=z_j$,
$L$~is the slope of the line segment, and $q_2(\alpha(t_m))$, $m=k\ldots,i$, are
approximations of $q_2$ at each $\alpha(t_m)$
obtained again from the interpolation of $q_2(z_r)$, $r=1,\ldots,M$, by a cubic~spline.
Note, $L = \frac{z_j - z_l}{t_i - t_k} > 0$ as $z_j > z_l, t_i > t_k$.
The energy for $\gamma$, which we denote by $E_0(\gamma)$, is then defined as the sum of
the energies over the linear components of the graph of~$\gamma$ with $\alpha$ in
\eqref{E:energy6} coinciding with $\gamma$ on each component. Thus, the problem of
minimizing \eqref{E:energy-discr} is then replaced by the approximately equivalent and
easier to solve problem of minimizing $E_0(\gamma)$ for~$\gamma\in\Gamma$.
\par For the purpose of efficiently computing $\gamma^*\in\Gamma$ of approximately
minimum energy, i.e., $\gamma^*\in\Gamma$ such that $\gamma=\gamma^*$ approximately
minimizes $E_0(\gamma)$ for $\gamma\in\Gamma$, we present an algorithm below that uses DP
on sets of grid points around graphs of estimates of $\gamma^*$, one set at a time,
each set containing the corner grid points $(1,1)$ and~$(N,M)$, each set defined each time
a new estimate of $\gamma^*$ is identified, each set coarser than the previous one, each set
contained in and associated with an essentially thin region, i.e., a strip, in the unit square
that contains the current estimate of~$\gamma^*$. In this algorithm, a general DP procedure,
Procedure~\emph{DP}, whose outline follows, is executed, for each strip as just described,
on the set $R$ of grid points associated with the strip, for the purpose of minimizing
$E_0(\gamma)$ (adjusted for $R$, see below) for $\gamma\in\Gamma$, $\gamma$ satisfying that
all of its vertices are in~$R$.
For such sets $R$ the computational cost is low (the search space is relatively small), and
their selection is such that it is highly likely the final DP solution is of approximately minimum
energy and therefore can be assumed to be the desired $\gamma^*$. Since the collection of
strips has the appearance of one single strip whose shape evolves as it mimics the shapes of
graphs of estimates of $\gamma^*$, we think of the collection as indeed being one single strip,
a dynamic strip that we call \emph{adapting} strip accordingly.
%In \cite{dogan-bernal-hagwood}, Dogan, Bernal and Hagwood proposed using DP on $R$
%in a strip of linear ($O(N)$) width around the diagonal of $[0,1]^2$ connecting planar
%points $(0,0)$ and $(1,1)$, for a fast DP algorithm.
Thus, we propose using an adapting strip as just described with a width that is constant
($O(1)$) as it evolves around graphs of estimates of $\gamma^*$. Obviously we do not know
$\gamma^*$, but can estimate it using DP solutions as the sets $R$ associated with the strips
 become coarser. However, before going into the specifics of our proposed algorithm, we will
describe Procedure~\emph{DP} operating on a generic $R$.
\par The set $R$ of labeled grid points can be any subset of the interior grid points plus
the corner grid points $(1,1)$, $(N,M)$. Given any such $R$, we denote by $\Gamma(R)$
the set of elements of $\Gamma$ with all vertices in $R$ (note, $R$ can have as few as three
points in which case $\Gamma(R)$ has only one element). Accordingly, with the energy
in \eqref{E:energy6} adjusted for~$R$ (see below) and $E_0(\gamma)$ for
$\gamma\in\Gamma(R)$ adjusted accordingly, given a positive integer $layrs$
(e.g., $layrs~=~5$) which determines the size of certain neighborhoods to be searched
(see below), then, based on DP, Procedure~\emph{DP} that follows, in $O(|R|)$ time,
will often (depending on $layrs$) compute $\gamma^*\in\Gamma(R)$ such that
$\gamma=\gamma^*$ approximately minimizes $E_0(\gamma)$ (adjusted for $R$, see below)
for~$\gamma\in\Gamma(R)$, $|R|$ the cardinality~of~$R$.
\par As the DP procedure progresses over the indices $(i,j)$ in $R$, it examines
function values on indices $(k,l)$ in a trailing neighborhood $N(i,j)$
of $(i,j)$.
%(see Figure~\ref{F:DP-strip} for a particular $R$ described below).
In the full DP, we would be examining all $(k,l)$ in $R$, $1 \leqslant k < i,
1 \leqslant l < j$. This has high computational cost, and is not necessary
for our applications. Using a much smaller square neighborhood $N(i,j)$ of
$\omega$ points ($\omega = layrs$) per side gives satisfactory results.
Thus, for each $(i,j)$
in $R$, we examine at most $\omega^2$ indices $(k,l)$ in the trailing neighborhood $N(i,j)$
(defined below). Then the overall time complexity is $O(\omega^2|R|)$.
We formally define $N(i,j)$~by
\begin{align*}
N(i,j) = \{
&(k,l) \in R:
k \ \mathrm{is\ one\ of\ } \omega \ \mathrm{largest\ indices\ } < i \\
&\mathrm{and} \
l \ \mathrm{is\ one\ of\ } \omega \ \mathrm{largest\ indices\ } < j \}.
\end{align*}
Note that in the case in which $N(i,j)$ happens to be empty then a grid point $(k,l)$
in $R$, $k<i$, $l<j$, perhaps $(k,l)=(1,1)$, is identified and $N(i,j)$ is
set to~$\{(k,l)\}$
\par The DP procedure follows. First, however, we clarify some implicit conventions
in the procedure logic. Grid points $(i,j)$ in $R$ are processed one at a time.
However, the main loop in the DP procedure takes place over the single index~$i$.
We process index $i$ in increasing order, and for each~$i$, each grid point $(i,j)$ in~$R$
is processed before moving to the next~$i$. Also in the procedure, pairs of indices $m_1$,
$m_2$ are retrieved from an index set~$\mathcal{M}$, satisfying $m_1<m_2$ with no other
index in~$\mathcal{M}$ greater than $m_1$ and less than $m_2$. This has the effect of
adjusting for $R$ energies computed according to~\eqref{E:energy6}.
\begin{tabbing}
12\=45\=78\=01\=34\=67\=90\=23\=56\=89\= \kill
{\small\bf procedure} \emph{DP}\\
\> $ E(1,1)=0$\\
\>{\small\bf for} each $(i,j)\not= (1,1)$ in $R$ in increasing order of $i$ {\small\bf do}\\
\>\>{\small\bf for} each $(k,l) \in N(i,j)$ {\small\bf do}\\
\>\>\> $\alpha =$ linear function, $\alpha(t_k)=z_l$, $\alpha(t_i)=z_j$\\
\>\>\> $L =$ slope of line segment~$\overline{(k,l)(i,j)}$\\
\>\>\> $\mathcal{M} = \{m:k\leq m\leq i, \exists (m,n) \in R \}$\\
\>\>\> $F_m= F(t_m,\alpha(t_m),L)$ for each $m\in \mathcal{M}$\\
\>\>\> $E_{(k,l)}^{(i,j)} = \frac{1}{2} \sum_{m_1,m_2 \in \mathcal{M}} (t_{m_2}-t_{m_1}) (F_{m_2} + F_{m_1})$\\
\>\>{\small\bf end for}\\
\>\> $E(i,j) = \min_{(k,l)\in N(i,j)} (E(k,l) + E_{(k,l)}^{(i,j)})$\\
\>\> $P(i,j) = \arg\min_{(k,l)\in N(i,j)} (E(k,l) + E_{(k,l)}^{(i,j)})$\\
\>{\small\bf end for}\\
{\small\bf end procedure}
\end{tabbing}
\par We note that $\gamma^*\in\Gamma(R)$ such that $\gamma=\gamma^*$ approximately minimizes
$E_0(\gamma)$ (adjusted for $R$) for $\gamma\in\Gamma(R)$ can then be obtained by backtracking
from $(N,M)$ to $(1,1)$ with pointer $P$ above. Accordingly,
Procedure~\emph{opt-diffeom} that follows, will produce $\gamma^*$ in the form
$\vec{\gamma}^* = (\gamma^*_{m})_{m=1}^N = (\gamma^*(t_m))_{m=1}^N$:
\begin{tabbing}
12\=45\=78\=01\=34\=67\=90\=23\=56\=89\= \kill
{\small\bf procedure} \emph{opt-diffeom}\\
\> $\gamma^*_{N} = 1$\\
\> $(i,j) = (N,M)$\\
\>{\small\bf while} $(i,j)\not= (1,1)$ {\small\bf do}\\
\>\> $(k,l) = P(i,j)$\\
\>\> $\gamma^*_{k}=z_l$\\
\>\>{\small\bf for} each integer $m, k<m<i$ {\small\bf do}\\
%\>\>\> $denom=t_i-t_k$\\
\>\>\> $\gamma^*_{m} = \frac{(t_i-t_m)}{(t_i-t_k)}z_l+
%\>\>\>\>\>\>
\frac{(t_m-t_k)}{(t_i-t_k)}z_j$\\
\>\>{\small\bf end}\\
\>\> $(i,j) = (k,l)$\\
\>{\small\bf end while}\\
{\small\bf end procedure}
\end{tabbing}
%
%\par In what follows, working with partitions (not necessarily uniform) of $[0,1]$,
%$\{t_l\}_{l=1}^{N}$, $\{z_l\}_{l=1}^{M}$, as previously described, we present
\par In what follows, we present the linear DP algorithm which we call \emph{adapt-DP},
based on DP restricted to an adapting strip, to compute an approximately optimal
diffeomorphism for the elastic registration of two curves in $d-$dimensional space. It
has parameters $layrs$, $lstrp$, set to small positive integers, say 5, 30, respectively.
Parameter $layrs$ is as previously described, while $lstrp$ is an additional
parameter that determines the width of the adapting strip (see below). Although
\emph{adapt-DP} is not guaranteed to be always successful (one or both of $layrs$,
$lstrp$ may be too small), it has been observed to produce convincing results for
every problem on which the algorithm has been applied, a few times after adjustments
to one or both parameters. The original ideas for this algorithm are described in
\cite{karypis,salvador} in the context of graph bisection and dynamic time~warping.
\par As presented below, \emph{adapt-DP} is essentially an iterative process that
restricts its search to the adapting strip around graphs of estimated solutions.
Each iteration culminates with the execution of Procedure~\emph{DP} for recursively
projecting a diffeomorphism obtained from a lower resolution grid to one of higher
resolution until full resolution is attained. For simplicity, we assume here
$N=M=2^n+1$ for some positive integer~$n$. Extending the algorithm to allow $N$,
$M$ to have any values is straightforward. Note, we don't assume partitions $\{t_l\}$,
$\{z_l\}$ are uniform. Finally, after the last execution of Procedure~\emph{DP} in
\emph{adapt-DP}, Procedure \emph{opt-diffeom} is performed to obtain, depending
on $layrs$ and $lstrp$, optimal $\gamma^*$ in~$\Gamma$.
Algorithm \emph{adapt-DP}~\mbox{follows:}
\begin{tabbing}
12\=45\=78\=01\=34\=67\=90\=23\=56\=89\= \kill
{\small\bf algorithm} \emph{adapt-DP} (DP algorithm)\\
2.\> $I(1) = J(1) = 1$\\
3.\> $P(N,M) = (1,1)$\\
\>{\small\bf for} $r=1$ {\small\bf to} $n$ {\small\bf do}\\
5.\>\> $NI=NJ=2^r+1$\\
6.\>\>{\small\bf for} $m=1$ {\small\bf to} $NI-1$ {\small\bf do}\\
7.\>\>\> $I(m+1) = m\cdot 2^{n-r}+1$\\
8.\>\>\> $r'_m = \frac{1}{2}(t_{I(m)} + t_{I(m+1)})$\\
\>\>{\small\bf end for}\\
\>\>{\small\bf for} $m=1$ {\small\bf to} $NJ-1$ {\small\bf do}\\
\>\>\> $J(m+1) = m\cdot 2^{n-r}+1$\\
12.\>\>\> $s'_m = \frac{1}{2}(z_{J(m)} + z_{J(m+1)})$\\
\>\>{\small\bf end for}\\
14.\>\> $r'_1=s'_1=0$\\
15.\>\> $r'_{NI-1}=s'_{NJ-1}=1$\\
\>\> $(i,j) = (N,M)$\\
\>\> $D=\emptyset$\\
18.\>\>{\small\bf while} $(i,j)\not= (1,1)$ {\small\bf do}\\
\>\>\> $(k,l) = P(i,j)$\\
**********************************************\\
20.\ Here below, for integers $m'$, $n'$, $1<m'<NI$,\\
21.\ $1<n'<NJ$, bin $B(m',n') \equiv$\\
22.\ $\{(x,y): r'_{m'-1}\leq x\leq r'_{m'}, s'_{n'-1}\leq y\leq s'_{n'}\}$\\
**********************************************\\
\>\>\> {\small\bf identify} bins $B(m',n')$, $1<m'<NI$,\\
\>\>\>$1<n'<NJ$, the interiors of which are\\
\>\>\>intersected by line segment $\overline{(i,j)(k,l)}$\\
\>\>\> $D'= \{(m',n'):\ \overline{(i,j)(k,l)}\cap B(m',n')\not=\emptyset\}$\\
\>\>\> $D=D\cup D'$\\
%\>\>\>$B(m',n')\not=\emptyset\}$\\
\>\>\> $(i,j) = (k,l)$\\
\>\>{\small\bf end while}\\
\>\> $R= \{(1,1),(N,M)\}$\\
31.\>\>{\small\bf for} each $(m',n')$ in $D$ {\small\bf do}\\
\>\>\> $i_0 = \max\{ 2,m'-lstrp\}$\\
\>\>\> $j_0 = \max\{ 2,n'-lstrp\}$\\
\>\>\> $R_1=\{(i,j): i=I(i'), i_0\leq i'\leq m', j=J(n')\}$\\
\>\>\> $R_2=\{(i,j): j=J(j'), j_0\leq j'\leq n', i=I(m')\}$\\
\>\>\> $R=R\cup R_1\cup R_2$\\
\>\>{\small\bf end for}\\
38.\>\>{\small\bf execute} procedure \emph{DP} on $R$\\
\>{\small\bf end for}\\
\>{\small\bf execute} procedure \emph{opt-diffeom} to obtain $\gamma^*$\\
{\small\bf end algorithm}
\end{tabbing}
In the outline of \emph{adapt-DP} above, we note in line~5, $NI$ starts equal
to $3$ (for $r=1$) and then it is essentially doubled at each iteration $r>1$
until it becomes equal to $N$ at the $n^{th}$ iteration. We note in line~2
and in line~7 inside the {\small\bf for} loop at line~6, the range of $I$ starts with 3
integers (for $r=1$) and then essentially doubles in size at each iteration
$r>1$, contains the previous range of~$I$ from preceding iteration, and is evenly
spread in the set $\{1,2,\ldots,N\}$ until it becomes this set. We note as well
from the well-known sum of a geometric series that since $N=2^n+1$ then the sum of the
$NI$'s, i.e., $(2^1+1)+(2^2+1)+\ldots+(2^n+1)$, is $O(N)$. Clearly, all of the above
applies to $NJ$, $M$, and the range of~$J$.
\begin{figure}
\centering
\begin{tabular}{ccc}
\includegraphics[width=0.25\textwidth]{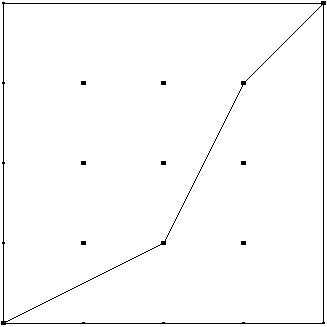}
&
\includegraphics[width=0.25\textwidth]{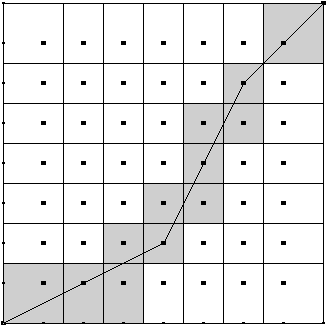}
&
\includegraphics[width=0.25\textwidth]{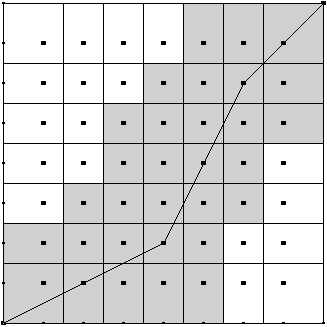}
\end{tabular}
\caption{\label{F:DP-wind}
On left is $\gamma^*$ from $2^{nd}$ iteration, $NI=NJ=2^2+1=5$.
In center, during $3^{rd}$ iteration, $NI=NJ=2^3+1=9$; shaded bins
are bins the interior of which $\gamma^*$ intersects.
On right, shaded bins form adapting strip in which next $\gamma^*$
is computed. Each shaded bin is within 2 bins ($lstrp=2$) of a bin
which is above, to the right of, or equal~to~it, and
whose interior has nonempty intersection with the current~$\gamma^*$.
}
\end{figure}
\par We note that the {\small\bf while} loop at line~18 identifies certain cells in the
Voronoi diagram \cite{voronoi} of the set of grid points $R'\equiv \{(i,j):i=I(m'),j=J(n')$,
$1<m'<NI$, $1<n'<NJ\}$ restricted to the unit square. Indeed bin $B(m',n')$ as defined
in lines~20-22, in terms of the computations in lines 8, 12, 14, 15, is exactly
the Voronoi cell of~$(I(m'),J(n'))$, and all such cells together partition the unit
square. Accordingly, with $\gamma^*$ encoded in~$P$ in line~3 ($r=1$) or in
line~38 ($r>1$) through the execution of Procedure~\emph{DP} in the previous iteration
($r-1$), it must be that every point in the graph of $\gamma^*$ is in some
bin~$B(m',n')$.  Thus, it then seems reasonable to say that a reliable region of
influence of $\gamma^*$ is the region around the graph of $\gamma^*$ formed by the union
of bins within a constant number of bins from the graph. Accordingly,
to be precise, a bin $B$ is part of this region if and only if there is a bin $B'$,
the interior of which the graph of~$\gamma^*$ intersects, $B$ within a constant
number ($lstrp$) of bins from $B'$, $B$ directly below or to the left of $B'$,
or $B$ equal to $B'$ (see Figure~\ref{F:DP-wind}). We note that identifying this
region is essentially accomplished in the {\small\bf while} loop at line~18 and the
{\small\bf for} loop at line~31, with the region understood to be the union of
bins or Voronoi cells $B(m',n')$ of grid points in $R$ at the end of the
{\small\bf for}~loop. Clearly, the region contains the graph of~$\gamma^*$,
and has the appearance of a strip whose shape evolves from one iteration to the next
as it closely mimics the shape of the graph of $\gamma^*$ (see Figure~\ref{F:DP-wind}),
and thus it is referred to as an adapting strip. Finally, we note that at the end of the
{\small\bf for}~loop, $\gamma^*$ in~$\Gamma(R)\subseteq \Gamma(R')$ encoded in~$P$
for the current iteration, is obtained in line~38 with Procedure~\emph{DP} restricted to the
region or adapting strip, a region that as just described depends on all
previous~$\gamma^*$ functions from previous iterations. The last $\gamma^*$
obtained is then, depending on $layrs$, optimal in $\Gamma(R)$, and,
depending on $layrs$ and $lstrp$, in~$\Gamma(R')$.
\par With $\gamma^*$ as above during the execution of the {\small\bf while} loop at line~18 for
iteration~$r$, we note that since $\gamma^*$ is in $\Gamma(R)$ then the number of bins
$B(m',n')$ whose interiors the graph of $\gamma^*$ intersects must be $O(NI+NJ)$, which
is also the time required to find them one linear component of the graph at a time.
Since $|R|$ at end of the {\small\bf for}~loop at line~31 is then
\mbox{$O(lstrp)\cdot O(NI+NJ)$}, i.e., $O(NI+NJ)$, the complexity of Procedure~\emph{DP}
at line~38 is then $O(NI+NJ)$, and since as mentioned above the sum of the $NI$'s
and $NJ$'s is $O(N)$ and $O(M)$, respectively, then the complexity of \emph{adapt-DP}
must be~$O(N+M)$, implying \emph{adapt-DP} is linear.
\section{\large Computation of Rotation for Rigid Alignment of Curves in
$d-$dimensional Space}
In this section, we describe and justify in a purely algebraic manner the usual algorithm for
computing an optimal rotation for the rigid alignment of two simple curves in $d-$dimensional space,
$d$ a positive integer, the curves not necessarily closed, each curve with a fixed starting point.
Here the computation of the alignment is based only on rotations (with $d\times d$ rotation matrices)
of one of the curves. We note, if $d$ equals~1, the only $1\times 1$ rotation matrix possible
is the one whose sole entry is~1. For simplicity we say that this matrix equals~1. We also note,
results presented here for the justification already appear in~\cite{lawrence},
although not in the context of shape~analysis.
\par We assume we have two simple curves in $d-$dimensional space represented by
functions $\beta_n: [0,1]\rightarrow \mathbb{R}^d$, $n=1,2$, that are absolutely continuous
and of unit length, and square-integrable functions $q_n: [0,1]\rightarrow \mathbb{R}^d$,
$\|q_n\|_{L^2}^2 = \int_0^1 \|q_n(t)\|^2 dt = 1$, $n=1,2$, the shape functions or SRVF's of
$\beta_n$, $n=1,2$, respectively, where $\|\cdot\|$ is the $d-$dimensional Euclidean norm
and $\|\cdot\|_{L^2}$ is the $L^2$~norm for functions in $L^2([0,1],\mathbb{R}^d)$.
Given $n$, $n=1,2$, we also assume $\beta_n(0)$ is the fixed starting point of~$\beta_n$.
Note, given an absolutely continuous function $\beta: [0,1]\rightarrow \mathbb{R}^d$ and
$R$ a rotation matrix in $\mathbb{R}^d$, if $q$ is the SRVF of $\beta$, then $Rq$ is the
SRVF of $R\beta$ (see~\cite{srivastava}).
Ideally, an optimal rotation matrix $R$ in $\mathbb{R}^d$ is found, i.e., a $d\times d$
orthogonal matrix $R$ with det$(R)=1$ (of determinant~1), that minimizes
%$
%R(\theta) = \left(
%\begin{smallmatrix}
%\cos(\theta) & \sin(\theta) \\
%-\sin(\theta) & \cos(\theta)\\
%\end{smallmatrix} \right)
%$
%
\begin{equation}\label{E:energy0}
E(R) = \int_0^1 \| q_1(t) - R q_2(t) \|^2 dt.
\end{equation}
For this purpose $E(R)$ in~\eqref{E:energy0} is rewritten as follows:
\begin{equation*}
E(R) = \|q_1\|_{L^2}^2 + \|q_2\|_{L^2}^2 - 2\int_0^1 q_1^T(t)R q_2(t)dt =
2 - 2\int_0^1 q_1^T(t)R q_2(t)dt
\end{equation*}
(note, $\|Rq_2(t)\|=\|q_2(t)\|$
since~$R$ is a rotation matrix). Then minimizing~\eqref{E:energy0} over all rotation
matrices~$R$ is equivalent to maximizing
\begin{equation*}
\int_0^1 q_1^T(t)R q_2(t)dt=\mathrm{tr}(RA^T),
\end{equation*}
where tr$(RA^T)$ is the trace of $RA^T$ and $A$ is the $d\times d$ matrix with entries
$A_{kj}=\int_0^1 q_{1k}(t)q_{2j}(t)dt$, where $q_{1k}(t)$, $q_{2j}(t)$ are
the $k^{th}$ and $j^{th}$ coordinates of $q_1(t)$ and $q_2(t)$, respectively, for each pair
$k,j=1,\ldots,d$. We refer to this matrix as the matrix~$A$ associated with~\eqref{E:energy0}.
\par As noted in Section~$2$, in practice, we need to solve a discretized version of the
problem. Thus, $\beta_1$ and $\beta_2$ are given as finite lists of points, say $N$ points per
curve for some integer $N>0$. For some partition $\{t_l\}_{l=1}^N$, $t_1=0<t_2<\ldots<t_N=1$,
(not necessarily uniform) of $[0,1]$, then for $n=1,2$, $\beta_n$ is given as $\beta_n(t_l)$,
$l=1,\ldots,N$. Similarly for $q_1$, $q_2$, except that for $l=1,\ldots,N$, $q_1(t_l)$ and
$q_2(t_l)$ are computed as described in Section~2. That the lists for $\beta_1$ and $\beta_2$,
and therefore for $q_1$ and $q_2$, must have the same number of points with the same partition
is dictated by the way optimal rotation matrices are computed as described below.
In what follows, for $l=1,\ldots,N$, $k=1,\ldots,d$, $j=1,\ldots,d$, $q_1^l$ is $q_1(t_l)$,
$q_2^l$ is $q_2(t_l)$, $q_{1k}^l$ is the $k^{th}$ coordinate of $q_1^l$, and $q_{2j}^l$
is the $j^{th}$ coordinate of $q_2^l$. We then discretize integral~\eqref{E:energy0} using
the trapezoidal rule:
\[ E^{discr}(R) = \sum_{l=1}^{N-1}1/2\,(t_{l+1}-t_l)(\|q_1(t_l)-Rq_2(t_l)\|^2
+\|q_1(t_{l+1})-Rq_2(t_{l+1})\|^2) \]
\begin{equation}\label{E:disc-energy0}
=\sum_{l=1}^{N}h_l\,\| q_1(t_l) - R q_2(t_l) \|^2
=\sum_{l=1}^{N}h_l\,\| q_1^l - R q_2^l \|^2,
\end{equation}
where $h_1=(t_2-t_1)/2$, $h_N=(t_N-t_{N-1})/2$, and for \mbox{$l=2,\ldots,N-1$},
$h_l=(t_{l+1}-t_{l-1})/2$.
Note, $h_l>0$ for each $l$, $l=1,\ldots,N$, and
$\sum_{l=1}^Nh_l=1$.\\ \smallskip\\
%
%With $h=1/(N-1)$, note that if the partition $\{t_l\}_{l=1}^N$ is uniform,
%then \eqref{E:disc-energy0} reduces to
%\begin{equation*}
%E^{discr}(R)=\sum_{l=1}^{N-1}h\,\| q_1(t_l) - R q_2(t_l) \|^2
%=\sum_{l=1}^{N-1}h\,\| q_1^l - R q_2^l \|^2.
%\end{equation*}
%
Note that minimizing \eqref{E:disc-energy0} over all rotation matrices~$R$ is an
instance of solving the so-called Wahba's problem \cite{wahba,markley} which is the problem
of minimizing $\sum_{l=1}^{N}w_l\,\| q_1^l - R q_2^l \|^2$ over all rotation
matrices~$R$, where the $w_l$'s are nonnegative weights.\\ \smallskip\\
Noting $\|Rq_2^l\| = \|q_2^l\|$, $l=1,\ldots,N$, we can rewrite
\eqref{E:disc-energy0} as follows
\begin{equation*}
E^{discr}(R) = \sum_{l=1}^{N}h_l\,(\|q_1^l\|^2 + \|q_2^l\|^2)
- 2\,\sum_{l=1}^{N}h_l\,((q_1^l)^T R q_2^l),
\end{equation*}
so that minimizing \eqref{E:disc-energy0} over all rotation matrices~$R$ is equivalent
to maximizing
\begin{equation}\label{E:trace}
\sum_{l=1}^{N}h_l\,(q_1^l)^T R q_2^l = \mathrm{tr}(RA^T),
\end{equation}
where $A$ is the $d\times d$ matrix with entries
$A_{kj} = \sum_{l=1}^{N}h_l\,q_{1k}^lq_{2j}^l$, for each pair $k,j = 1,\ldots,d$.
We refer to this matrix as the matrix~$A$ associated with~\eqref{E:disc-energy0}.
Focusing our attention on this matrix, as opposed to doing it on the matrix~$A$
associated with~\eqref{E:energy0}, then an optimal rotation matrix $R$
for~\eqref{E:trace}, thus for~\eqref{E:disc-energy0}, can be computed from the
singular value decomposition of $A$ or, more precisely, with the Kabsch-Umeyama
algorithm \cite{kabsch1,kabsch2,umeyama}
(see Algorithm Kabsch-Umeyama below, where $\mathrm{diag}\{s_1,\ldots,s_d\}$
is the $d\times d$ diagonal matrix with numbers $s_1,\ldots,s_d$ as the elements
of the diagonal, in that order running from the upper left to the lower right of
the matrix). A {\em singular value decomposition} (SVD)~\cite{lay} of~$A$ is a
representation of the form $A=USV^T$, where $U$ and $V$ are~$d\times d$ orthogonal
matrices and $S$ is a $d\times d$ diagonal matrix with the singular values of $A$, which
are nonnegative real numbers, appearing in the diagonal of $S$ in descending order,
from  the upper left to the lower right of~$S$. Note that any matrix, not
necessarily square, has a singular value decomposition, not necessarily unique~\cite{lay}.
Note as well that an optimal rotation matrix $R$ for integral~\eqref{E:energy0} can also be
computed in a similar manner using the matrix $A$ associated with~\eqref{E:energy0}.
However, in what follows, for the obvious reasons, we focus our attention on
the discretized integral~\eqref{E:disc-energy0}. Finally, note that applications of
the Kabsch-Umeyama algorithm similar to the one just decribed here can be found
in~\mbox{\cite{dogan,srivastava}}.
%\begin{algorithm}
%\caption{Computing optimal $R$}
%\label{A:kabsch}
\begin{algorithmic}
\STATE \noindent\rule{13cm}{0.4pt}
\STATE {\bf Algorithm Kabsch-Umeyama} (KU algorithm)
\STATE \noindent\rule[.1in]{13cm}{0.4pt}
\STATE Set $h_1=(t_2-t_1)/2$, $h_N=(t_N-t_{N-1})/2$, $h_l=(t_{l+1}-t_{l-1})/2$ for\\ $l=2,\ldots,N-1$.
\STATE Set $q^l_{1k}$ equal to $k^{th}$ coordinate of $q_1(t_l)$ for $l=1,\ldots,N$, $k=1,\ldots,d$.
\STATE Set $q^l_{2j}$ equal to $j^{th}$ coordinate of $q_2(t_l)$ for $l=1,\ldots,N$, $j=1,\ldots,d$.
\STATE Compute $A_{kj} = \sum_{l=1}^{N}h_l\,q_{1k}^l q_{2j}^l$ for each pair $k,j=1,\ldots,d$.
\STATE Identify $d\times d$ matrix $A$ with entries $A_{kj}$ for each pair $k,j=1,\ldots,d$.
%\STATE {\small\bf if} $d=1$ {\small\bf then} set $R=1$.
%\STATE {\small\bf else}
\STATE Compute SVD of $A$, i.e., identify $d\times d$ matrices $U$, $S$, $V$, so that
\STATE $A = U S V^T$ in the SVD sense.
\STATE Set $s_1= \ldots = s_{d-1}=1$.
\STATE {\small\bf if} $\det(UV) > 0$ {\small\bf then} set $s_d=1$.
\STATE {\small\bf else} set $s_d=-1$. {\small\bf end if}
\STATE Set $\tilde{S} = \mathrm{diag}\{s_1,\ldots,s_d\}$.
\STATE Compute and return $R = U \tilde{S} V^T$ and $maxtrace=\mathrm{tr}(RA^T)$.
%\STATE {\small\bf end if}
\STATE \noindent\rule{13cm}{0.4pt}
\end{algorithmic}
%\end{algorithm}
\smallskip
\par We note that if $d=1$, the KU algorithm still computes $R$ and
$maxtrace$, with the resulting $R$ always equal to~1.
We should also note that alternative methods for $d=2,\ 3$ to the Kabsch-Umeyama
algorithm, i.e., that do not use the SVD to compute a rotation matrix $R$ that maximizes
$\mathrm{tr}(RM)$ over all $d\times d$ rotations matrices, $M$ a $d\times d$ matrix,
have been presented in~\cite{bernal3}.
%\smallskip
\\ \smallskip
\par Now we justify the KU algorithm using exclusively simple concepts from linear algebra,
mostly in the proof of the proposition that follows. The algorithm has been previously justified in
\cite{kabsch1,kabsch2,umeyama}, however the justifications there are not totally algebraic as they
are based on the optimization technique of Langrange multipliers. The justification here already
appears in~\cite{lawrence}, however it was not developed there in the context of shape analysis.
Accordingly we develop it here with that context in mind (the matrix A in the outline of the
KU algorithm above is defined in that context) but mainly for self-containment and clarity.
\smallskip\\
{\bf Proposition 1:} If $D=\mathrm{diag}\{\sigma_1,\ldots,\sigma_d\}$, $\sigma_j\geq 0$,
$j=1,\ldots,d$, and $W$ is a $d\times d$ orthogonal matrix, then\\
1. tr$(WD)\leq\sum_{j=1}^d \sigma_j$.\\
2. If $B$ is a $d\times d$ orthogonal matrix, $S=B^TDB$, then $\mathrm{tr}(WS)\leq \mathrm{tr}(S)$.\\
3. If det$(W)=-1$, $\sigma_d\leq \sigma_j$, $j=1,\ldots,d-1$, then
$\mathrm{tr}(WD)\leq\sum_{j=1}^{d-1}\sigma_j-\sigma_d$.
\smallskip\\
{\bf Proof:} Since $W$ is orthogonal and if $W_{kj}$, $k,j=1,\ldots,d$, are the
entries of $W$, then, in particular, $W_{jj}\leq 1$, $j=1,\ldots,d$, so that\\
$\mathrm{tr}(WD)=\sum_{j=1}^d W_{jj}\sigma_j\leq \sum_{j=1}^d \sigma_j$, and therefore 1. holds.\\
Accordingly, assumming $B$ is a $d\times d$ orthogonal matrix, since $BWB^T$ is also
orthogonal, it follows from 1. that\\
$\mathrm{tr}(WS)=\mathrm{tr}(WB^TDB)=\mathrm{tr}(BWB^TD)\leq\sum_{j=1}^d\sigma_j=\mathrm{tr}(D)
=\mathrm{tr}(S)$, and therefore 2. holds.\\
If det$(W)=-1$, we show next that a $d\times d$ orthogonal matrix $B$ can be
identified so that with $\bar{W}=B^TWB$, then
$\bar{W}= \left( \begin{smallmatrix}
W_0 & O\\ O^T & -1\\
\end{smallmatrix} \right)$,
$W_0$ interpreted as the upper leftmost $d-1\times d-1$ entries of $\bar{W}$ and as a
$d-1\times d-1$ matrix as well; $O$ interpreted as a vertical column or vector of $d-1$ zeroes.\\
With $I$ as the $d\times d$ identity matrix, then det$(W)=-1$ implies
$\mathrm{det}(W+I)=-\mathrm{det}(W)\mathrm{det}(W+I)=-\mathrm{det}(W^T)\mathrm{det}(W+I)=
-\mathrm{det}(I+W^T)=-\mathrm{det}(I+W)$ which implies det$(W+I)=0$
so that $x\not=0$ exists in $\mathbb{R}^d$ with $Wx=-x$. It also follows then that
$W^TWx=W^T(-x)$ which gives $x=-W^Tx$ so that $W^Tx=-x$ as well.\\
Letting $b_d=x$, vectors $b_1,\ldots,b_{d-1}$ can be obtained so that $b_1,\ldots,b_d$ form
a basis of~$\mathbb{R}^d$, and by the Gram-Schmidt process starting with $b_d$, we may
assume $b_1,\ldots,b_d$ form an orthonormal basis of $\mathbb{R}^d$ with $Wb_d=W^Tb_d=-b_d$.
Letting $B=(b_1,\ldots,b_d)$, interpreted as a $d\times d$ matrix with columns $b_1,\ldots,b_d$,
in that order, then it follows that $B$ is orthogonal, and with $\bar{W}=B^TWB$ and
$W_0$, $O$ as previously described, noting
$B^TWb_d=B^T(-b_d)= \left( \begin{smallmatrix} O\\ -1\\ \end{smallmatrix} \right)$ and
$b_d^TWB=(W^Tb_d)^TB=(-b_d)^TB=(O^T \, -1)$,
then $\bar{W}= \left( \begin{smallmatrix}
W_0 & O\\ O^T & -1\\
\end{smallmatrix} \right)$. Note, $\bar{W}$ is orthogonal
and therefore so is the $d-1\times d-1$ matrix~$W_0$.\\
Let $S=B^TDB$ and write
$S = \left( \begin{smallmatrix}
S_0 & a\\ b^T & \gamma\\
\end{smallmatrix} \right)$,
$S_0$ interpreted as the upper leftmost $d-1\times d-1$ entries of $S$ and as a
$d-1\times d-1$ matrix as well; $a$ and $b$ interpreted as vertical columns or vectors of
$d-1$ entries, and $\gamma$ as a scalar.\\
Note, $\mathrm{tr}(WD)=\mathrm{tr}(B^TWDB)=\mathrm{tr}(B^TWBB^TDB)=\mathrm{tr}(\bar{W}S)$,
so that $\bar{W}S=$
$\left( \begin{smallmatrix}
W_0 & O\\ O^T & -1\\
\end{smallmatrix} \right)$
$\left( \begin{smallmatrix}
S_0 & a\\ b^T & \gamma\\
\end{smallmatrix} \right)=$
$\left( \begin{smallmatrix}
W_0S_0 & W_0a\\ -b^T & -\gamma\\
\end{smallmatrix} \right)$
gives $\mathrm{tr}(WD)=\mathrm{tr}(W_0S_0)-\gamma$.\\
We show $\mathrm{tr}(W_0S_0)\leq\mathrm{tr}(S_0)$. For this purpose let
$\hat{W}= \left( \begin{smallmatrix}
W_0 & O\\ O^T & 1\\
\end{smallmatrix} \right)$,
$W_0$ and $O$ as above. Since $W_0$ is orthogonal, then clearly $\hat{W}$ is a $d\times d$
orthogonal matrix, and by 2. $\mathrm{tr}(\hat{W}S)\leq \mathrm{tr}(S)$
so that $\hat{W}S=$
$\left( \begin{smallmatrix}
W_0 & O\\ O^T & 1\\
\end{smallmatrix} \right)$
$\left( \begin{smallmatrix}
S_0 & a\\ b^T & \gamma\\
\end{smallmatrix} \right)=$
$\left( \begin{smallmatrix}
W_0S_0 & W_0a\\ b^T & \gamma\\
\end{smallmatrix} \right)$
gives $\mathrm{tr}(W_0S_0)+\gamma=\mathrm{tr}(\hat{W}S)\leq\mathrm{tr}(S)=\mathrm{tr}(S_0)+\gamma$.
Thus, $\mathrm{tr}(W_0S_0)\leq\mathrm{tr}(S_0)$.\\
Note, $\mathrm{tr}(S_0)+\gamma = \mathrm{tr}(S)= \mathrm{tr}(D)$, and if $B_{kj}$,
$k,j=1,\ldots,d$ are the entries of $B$, then $\gamma = \sum_{k=1}^d B_{kd}^2\sigma_k$,
a convex combination of the $\sigma_k$'s, so that $\gamma\geq\sigma_d$.
It then follows that\\ $\mathrm{tr}(WD)=\mathrm{tr}(W_0S_0)-\gamma\leq\mathrm{tr}(S_0)-\gamma=
\mathrm{tr}(D)-\gamma-\gamma\leq\sum_{j=1}^{d-1}\sigma_j-\sigma_d$, and therefore
3. holds. \ $\Box$
%\\ \smallskip\par
\smallskip\\
Finally, the following theorem, a consequence of Proposition~1, justifies the
Kabsch-Umeyama algorithm.
\smallskip\\
{\bf Theorem:} Given a $d\times d$ matrix $M$, let $U$, $S$, $V$ be $d\times d$ matrices such that
the singular value decomposition of $M$ gives $M=USV^T$. If det$(UV^T)>0$, then $R=UV^T$
maximizes tr$(RM^T)$ over all $d\times d$ rotation matrices~$R$. Otherwise, if det$(UV^T)<0$, with
$\tilde{S}=\mathrm{diag}\{s_1,\ldots,s_d\}$, $s_1 = \ldots = s_{d-1}=1,\ s_d=-1$,
then $R=U\tilde{S}V^T$ maximizes tr$(RM^T)$ over all $d\times d$ rotation matrices~$R$.
\smallskip\\
{\bf Proof:} Let $\sigma_j$, $j=1,\ldots,d$, $\sigma_1\geq \sigma_2\geq\ldots\geq\sigma_d\geq 0$,
be the singular values of~$M$, so that $S=\mathrm{diag}\{\sigma_1,\ldots,\sigma_d\}$.\\
Assume det$(UV^T)>0$. If $R$ is any rotation matrix, then $R$ is orthogonal.
From 1. of Proposition~1 since $U^TRV$ is orthogonal, then\\
$\mathrm{tr}(RM^T)=\mathrm{tr}(RVSU^T)=\mathrm{tr}(U^TRVS)\leq\sum_{j=1}^d\sigma_j.$\\
On the other hand, if $R=UV^T$, then $R$ is clearly orthogonal, \mbox{det$(R)=1$}, and
$\mathrm{tr}(RM^T)=\mathrm{tr}(UV^TVSU^T)=\mathrm{tr}(USU^T)=\mathrm{tr}(S)
=\sum_{j=1}^d\sigma_j.$\\
Thus, $R=UV^T$ maximizes tr$(RM^T)$ over all $d\times d$ rotation matrices~$R$.\\
Finally, assume det$(UV^T)<0$. If $R$ is any rotation matrix, then $R$ is orthogonal
and~det$(R)=1$. From 3. of Proposition~1 since $U^TRV$ is orthogonal and det$(U^TRV)=-1$, then\\
$\mathrm{tr}(RM^T)=\mathrm{tr}(RVSU^T)=\mathrm{tr}(U^TRVS)\leq\sum_{j=1}^{d-1}\sigma_j-\sigma_d.$\\
On the other hand, if $R=U\tilde{S}V^T$, then $R$ is clearly orthogonal, \mbox{det$(R)=1$}, and
$\mathrm{tr}(RM^T)=\mathrm{tr}(U\tilde{S}V^TVSU^T)=
\mathrm{tr}(U\tilde{S}SU^T)=\mathrm{tr}(\tilde{S}S)
=\sum_{j=1}^{d-1}\sigma_j-\sigma_d.$\\
Thus, $R=U\tilde{S}V^T$ maximizes tr$(RM^T)$ over all $d\times d$ rotation matrices~$R$.  \ $\Box$
\\ \smallskip
\par In the rest of this section, although not exactly related to the goal of this paper, for the
sake of completeness, we show how another problem of interest reduces to the problem just solved
above so that it can then be solved
with the Kabsch-Umeyama algorithm. The problem of interest is the so-called orientation-preserving
rigid motion problem. With $q_1$, $q_2$, $h_l$, $q_1^l$, $q_2^l$, $l=1,\ldots,N$, as above, the
problem is then that of finding an orientation-preserving rigid motion~$\phi$ of $\mathbb{R}^d$
that minimizes
\begin{equation}\label{E:rigid-energy0}
\Delta(\phi)=\sum_{l=1}^{N}h_l\,\| q_1^l - \phi(q_2^l) \|^2.
\end{equation}
An affine linear function $\phi$, $\phi:\mathbb{R}^d\rightarrow\mathbb{R}^d$, is a rigid
motion of $\mathbb{R}^d$ if it is of the form $\phi(q)= Rq + t$ for $q\in\mathbb{R}^d$, where $R$
is a $d\times d$ orthogonal matrix, and $t\in\mathbb{R}^d$. The rigid motion $\phi$ is orientation
preserving if det$(R)=1$.
We note that the justification of the reduction here, apparently already appears in~\cite{lawrence}.
However there the problem being reduced does not involve a partition of $[0,1]$ the way it does here,
a partition that can be either uniform or nonuniform, although, if so desired, it is not hard to show
that a partition can actually be associated with the problem in~\cite{lawrence}, a partition that must
be uniform. Thus the problem being reduced here is more general than the problem in~\cite{lawrence}.
\\ \smallskip\\
Let $\bar{q}_1$ and  $\bar{q}_2$ denote the centroids of the discretized $q_1$ and $q_2$,
respectively:
\begin{equation*}
\bar{q}_1=\sum_{l=1}^N h_l\,q_1^l
\ \ \ \ \mathrm{and}\ \ \ \ \bar{q}_2=\sum_{l=1}^N h_l\,q_2^l.
\end{equation*}
The following proposition shows, in particular, that $\phi(\bar{q}_2) = \bar{q}_1$ if $\phi$
minimizes~\eqref{E:rigid-energy0} over either the set of all rigid motions of $\mathbb{R}^d$
or the smaller set of rigid motions of $\mathbb{R}^d$ that are orientation preserving.
\smallskip\\
{\bf Proposition 2:} Let $\phi$ be a rigid motion of $\mathbb{R}^d$ with
$\phi(\bar{q}_2)\not=\bar{q}_1$ and
define an affine linear function $\tau$, $\tau: \mathbb{R}^d\rightarrow\mathbb{R}^d$,
$\tau(q) = \phi(q)-\phi(\bar{q}_2)+\bar{q}_1$ for $q\in\mathbb{R}^d$. Then~$\tau$~is a rigid
motion of~$\mathbb{R}^d$, $\tau(\bar{q}_2)=\bar{q}_1$, $\Delta(\tau) < \Delta(\phi)$, and if
$\phi$ is orientation preserving, then so is~$\tau$.
\smallskip\\
{\bf Proof:} Clearly $\tau(\bar{q}_2)=\bar{q}_1$.
Let $R$ be a $d\times d$ orthogonal matrix and $t\in\mathbb{R}^d$ such that
$\phi(q) = Rq+t$ for $q$ in~$\mathbb{R}^d$. Then $\tau(q)=Rq-R\bar{q}_2+\bar{q}_1$
so that $\tau$ is a rigid motion of $\mathbb{R}^d$, $\tau$ is orientation preserving if $\phi$ is,
and for each $l$, $l=1,\ldots,N$, we have
\par \,
\par \,
\hspace*{-.25in}$\|q_1^l-\phi(q_2^l)\|^2-\|q_1^l-\tau(q_2^l)\|^2\\
=(q_1^l-Rq_2^l-t)^T(q_1^l-Rq_2^l-t)$
\par $ -(q_1^l-Rq_2^l+R\bar{q}_2-\bar{q}_1)^T(q_1^l-Rq_2^l+R\bar{q}_2-\bar{q}_1)\\
=(q_1^l-Rq_2^l)^T(q_1^l-Rq_2^l)-2(q_1^l-Rq_2^l)^Tt+t^Tt-(q_1^l-Rq_2^l)^T(q_1^l-Rq_2^l)$
\par $-2(q_1^l-Rq_2^l)^T(R\bar{q}_2-\bar{q}_1)-(R\bar{q}_2-\bar{q}_1)^T(R\bar{q}_2-\bar{q}_1)\\
=-2(q_1^l-Rq_2^l)^Tt+t^Tt-2(q_1^l-Rq_2^l)^T(R\bar{q}_2-\bar{q}_1)$
\par $-(R\bar{q}_2-\bar{q}_1)^T(R\bar{q}_2-\bar{q}_1)
+2(R\bar{q}_2-\bar{q}_1)^Tt-2(R\bar{q}_2-\bar{q}_1)^Tt+t^Tt-t^Tt\\
=2(Rq_2^l-q_1^l)^Tt+2t^Tt+2(Rq_2^l-q_1^l)^T(R\bar{q}_2-\bar{q}_1)+ 2(R\bar{q}_2-\bar{q}_1)^Tt$
\par $-((R\bar{q}_2-\bar{q}_1)^T(R\bar{q}_2-\bar{q}_1)+2(R\bar{q}_2-\bar{q}_1)^Tt+t^Tt)\\
=2(Rq_2^l-q_1^l+t)^T(R\bar{q}_2-\bar{q}_1+t)-(R\bar{q}_2-\bar{q}_1+t)^T(R\bar{q}_2-\bar{q}_1+t)$.
\\ \smallskip\\
It then follows that
\par \,
\par \,
\hspace*{-.25in}$\Delta(\phi)-\Delta(\tau) = \sum_{l=1}^N h_l\|q_1^l-\phi(q_2^l)\|^2-
\sum_{l=1}^N h_l\|q_1^l-\tau(q_2^l)\|^2\\
= \sum_{l=1}^N h_l(\|q_1^l-\phi(q_2^l)\|^2-\|q_1^l-\tau(q_2^l)\|^2)\\
=\sum_{l=1}^N h_l(2(Rq_2^l-q_1^l+t)^T(R\bar{q}_2-\bar{q}_1+t)-
(R\bar{q}_2-\bar{q}_1+t)^T(R\bar{q}_2-\bar{q}_1+t))\\
= 2(R\sum_{l=1}^N h_l\,q_2^l-\sum_{l=1}^N h_l\,q_1^l
+t\,\sum_{l=1}^N h_l)^T(R\bar{q}_2-\bar{q}_1+t)$
\par $-(R\bar{q}_2-\bar{q}_1+t)^T(R\bar{q}_2-\bar{q}_1+t)\sum_{l=1}^N h_l\\
= 2(R\bar{q}_2-\bar{q}_1+t)^T(R\bar{q}_2-\bar{q}_1+t)-
(R\bar{q}_2-\bar{q}_1+t)^T(R\bar{q}_2-\bar{q}_1+t)\\
=(R\bar{q}_2-\bar{q}_1+t)^T(R\bar{q}_2-\bar{q}_1+t)=\|R\bar{q}_2-\bar{q}_1+t\|^2
=\|\phi(\bar{q}_2)-\bar{q}_1\|^2>0. \ \Box$
\\ \smallskip
\par Finally, the following corollary is a consequence of Proposition 2.
Here $p_1^l=q_1^l-\bar{q}_1$, $p_2^l=q_2^l-\bar{q}_2$, for $l=1,\ldots,N$,
and if $\bar{p}_1=\sum_{l=1}^Nh_l\,p_1^l$, $\bar{p}_2=\sum_{l=1}^Nh_l\,p_2^l$,
then clearly $\bar{p}_1=\bar{p}_2=0$.
It shows that the problem of finding an orientation-preserving rigid motion~$\phi$
of~$\mathbb{R}^d$ that minimizes~\eqref{E:rigid-energy0} can be reduced to the problem of
finding a $d\times d$ rotation matrix~$R$ that minimizes
$\sum_{l=1}^{N}h_l\,\| p_1^l - R p_2^l \|^2$ which, of course, then can be solved with the
Kabsch-Umeyama algorithm.
\smallskip\\
{\bf Corollary:} Let $\hat{R}$ be such that $R=\hat{R}$ minimizes
$\sum_{l=1}^Nh_l\,\|p_1^l-Rp_2^l\|^2$ over all $d\times d$ rotation matrices~$R$.
Let $\hat{t}=\bar{q}_1-\hat{R}\bar{q}_2$, and let $\hat{\phi}$ be given by
$\hat{\phi}(q)=\hat{R}q+\hat{t}$ for $q\in\mathbb{R}^d$. Then $\phi=\hat{\phi}$ minimizes
$\sum_{l=1}^Nh_l\,\|q_1^l-\phi(q_2^l)\|^2$ over all orientation-preserving rigid
motions $\phi$ of~$\mathbb{R}^d$.
\smallskip\\
{\bf Proof:} One such $\hat{R}$ can be computed with the Kabsch-Umeyama algorithm.
By Proposition~2, in order to minimize $\sum_{l=1}^Nh_l\,\|q_1^l-\phi(q_2^l)\|^2$ over all
orientation-preserving rigid motions $\phi$ of~$\mathbb{R}^d$, it suffices to do it over those
for which $\phi(\bar{q}_2)=\bar{q}_1$. Therefore, it suffices to minimize
$\sum_{l=1}^Nh_l\,\|q_1^l-(Rq_2^l+t)\|^2$ with $t=\bar{q}_1-R\bar{q}_2$ over all $d\times d$
rotation matrices~$R$, i.e., it suffices to minimize
$$\sum_{l=1}^Nh_l\|q_1^l-Rq_2^l-\bar{q}_1+R\bar{q}_2\|^2 =
\sum_{l=1}^Nh_l\,\|(q_1^l-\bar{q}_1)-R(q_2^l-\bar{q}_2)\|^2$$
over all $d\times d$ rotation matrices~$R$. But minimizing the last expression is equivalent
to minimizing $\sum_{l=1}^Nh_l\,\|p_1^l-Rp_2^l\|^2$ over all $d\times d$ rotation matrices~$R$.
Since $R=\hat{R}$ is a solution to this last problem, it then follows that
$R=\hat{R}$ minimizes $\sum_{l=1}^Nh_l\,\|q_1^l-Rq_2^l-\bar{q}_1+R\bar{q}_2\|^2$
$=\sum_{l=1}^Nh_l\,\|q_1^l -(Rq_2^l+t)\|^2$ with $t=\bar{q}_1-R\bar{q}_2$ over all
$d\times d$ rotation matrices~$R$. Consequently, if $\hat{t}=\bar{q}_1-\hat{R}\bar{q}_2$,
and $\hat{\phi}$ is given by $\hat{\phi}(q)=\hat{R}q+\hat{t}$ for $q\in\mathbb{R}^d$,
then $\phi=\hat{\phi}$ clearly minimizes $\sum_{l=1}^Nh_l\,\|q_1^l-\phi(q_2^l)\|^2$ over
all orientation-preserving rigid motions $\phi$ of~$\mathbb{R}^d$.~$\Box$
\section{\large Computation of the Elastic Shape Distance between Two Curves in
$d-$dimensional Space}
Let $\beta_1$, $\beta_2$, $q_1$, $q_2$ be as above, i.e.,
$\beta_n: [0,1]\rightarrow \mathbb{R}^d$, $n=1,2$, are absolutely continuous functions
representing simple curves in~$\mathbb{R}^d$ of unit length, and $q_n: [0,1]\rightarrow \mathbb{R}^d$,
$n=1,2$, are square-integrable functions that are the shape functions or SRVF's of $\beta_n$,
$n=1,2$, respectively. In the case one of the curves is closed, say $\beta_1$, then, in particular,
$q_1$ is interpreted to be a periodic function from $\mathbb{R}$ into $\mathbb{R}^d$ so that
$q_1(t+1)=q_1(t)$ for all vaues of~$t$. Define a finite subset $K$ of $[0,1]$ as follows.
If neither curve is closed let $K=\{0\}$. Otherwise, assume the curve represented by $\beta_1$,
the first curve, is closed, and for an integer $T>0$ choose numbers $\hat{t}_l$, $l=1,\ldots,T$,
in $[0,1]$, $\hat{t}_1=0<\ldots<\hat{t}_T=1$, so that $B=\{\beta_1(\hat{t}_l),l=1,\ldots,T\}$ is
reasonably dense in $\beta_1$ (by joining consecutive points in $B$ with line segments, the
resulting piecewise linear curve should be a reasonable approximation of~$\beta_1$).
Accordingly, let $K=\{\hat{t}_1,\ldots,\hat{t}_T\}$ in this case. Either way, $K$ is a finite
subset of~$[0,1]$ and $\{\beta_1(t), t\in K\}$ is interpreted to be the set of starting points
of~$\beta_1$.  With $SO(d)$ as the set of $d\times d$ rotation matrices,
and $\Gamma$ as the set of diffeomorphisms of $[0,1]$ onto itself ($\gamma(0)=0$, $\gamma(1)=1$,
$\dot{\gamma}>0$, for~$\gamma\in\Gamma$), given $t_0\in K$, $R\in SO(d)$, $\gamma\in\Gamma$,
a mismatch energy $E(t_0,R,\gamma)$ is defined by
\begin{equation}\label{E:dist-energy}
E(t_0,R,\gamma) = \int_0^1 \| \sqrt{\dot{\gamma}(t)} R q_1(t_0 + \gamma(t)) - q_2(t) \|^2 dt
\end{equation}
which as noted in \cite{dogan}, for the purpose of minimizing it with respect to $t_0$, $R$
and $\gamma$, without any loss of generality, can be reformulated as
\begin{equation}\label{E:dist2-energy}
E(t_0,R,\gamma) = \int_0^1 \| R q_1(t_0+t) - \sqrt{\dot{\gamma}(t)} q_2(\gamma(t)) \|^2 dt
\end{equation}
which allows for the second curve to be reparametrized while the first one is rotated.
Note, $SO(d)=\{1\}$ if $d=1$. \smallskip\\
As established in \cite{srivastava,srivastava2}, given $t_0\in K$, $R\in SO(d)$, $\gamma\in\Gamma$,
so that the triple $(t_0,R,\gamma)$ is a global minimizer of \eqref{E:dist-energy},
then $E(t_0,R,\gamma)$ can be interpreted to be the elastic shape distance between $\beta_1$
and $\beta_2$, and the elastic registration of $\beta_1$ and $\beta_2$ is obtained by
reparametrizing and rotating $\beta_1$ with $\gamma$ and $R$, respectively, with starting point
$\beta_1(t_0)$. Note, if $\beta_1$ is closed, we assume $\dot{\gamma}(0)=\dot{\gamma}(1)$ to ensure
the periodicity of $\sqrt{\dot{\gamma}(t)} R q_1(t_0 + \gamma(t))$ for~$t\in [0,1]$. Similarly, if
$(t_0,R,\gamma)$ is a global minimizer of
\eqref{E:dist2-energy} then again $E(t_0,R,\gamma)$ is interpreted to be the elastic shape distance
between $\beta_1$ and $\beta_2$, and the elastic registration of $\beta_1$ and $\beta_2$ is
obtained by rotating $\beta_1$ with $R$, with starting point $\beta_1(t_0)$, and reparametrizing
$\beta_2$ with $\gamma$.
\par As noted in Section 2 and Section 3, in practice, we work with curves $\beta_1$, $\beta_2$,
given as finite lists of points in the curves. Unless otherwise specified,
here and in what follows, for simplicity, although mostly out of necessity as will be made clear
below, we assume that $\beta_1$ and $\beta_2$ are discretized by the same partition of~$[0,1]$.
Even if at first they are not, this can be easily accomplished by interpolating with a cubic spline
one or both curves.
Accordingly, for some integer $N>0$, and a partition
of~$[0,1]$, $\{t_l\}_{l=1}^{N}$, $t_1=0<t_2<\ldots<t_N=1$, for $n=1,2$, the curve $\beta_n$ is given
as a list of $N$ points in the curve, where for $l=1,\ldots,N$, $\beta_n(t_l)$ is the $l^{th}$ point
in the list for~$\beta_n$. Similarly for $q_1$, $q_2$, except that for $l=1,\ldots,N$, $q_1(t_l)$
and $q_2(t_l)$ are computed as described in Section~2. That the lists for
$\beta_1$ and $\beta_2$, and therefore for $q_1$ and $q_2$, must have the same number of points with
the same partition is dictated by the way optimal rotation matrices are computed as described in
Section~3. If neither curve is closed so that each curve has a fixed starting point, then the partition
does not have to be uniform as pointed out in Section~3. The starting points of $\beta_1$ and $\beta_2$
are then $\beta_1(t_1)=\beta_1(0)$ and $\beta_2(t_1)=\beta_2(0)$, respectively, and as pointed out above
the finite subset $K$ of~$[0,1]$ defined above must equal~$\{0\}=\{t_1\}$. However if one curve is closed,
assumed to be $\beta_1$, then $\beta_2$ has a fixed starting point which is~$\beta_2(0)$, and any point
in $\beta_1$ can then be treated as a starting point of~$\beta_1$. In this case, for simplicity, $K$ is
purposely chosen to equal $\{t_1,\ldots,t_{N-1}\}$ or a subset of it, so that $B=\{\beta_1(t), t\in K\}$
can be assumed to be reasonably dense in~$\beta_1$ and therefore essentially of size~$O(N)$. $B$~is then
interpreted to be the set of starting points of the curve, and if the partition $\{t_l\}_{l=1}^N$ is
uniform, an optimal rotation matrix can then be easily computed as described in Section~3 for each
starting point in~$B$. We note that if $K$ equals $\{t_1,\ldots,t_{N-1}\}$, then the partition must indeed
be uniform. To see this we note that, in particular, for each $m$, $m=1,\ldots,N-1$, it should be that
$t_m+t_2=t_{m+1}$ so that $t_{m+1}-t_m=t_2=t_2-t_1$ and therefore the partition is uniform. Even if
$K$ is not all of $\{t_1,\ldots,t_{N-1}\}$ but a subset of it so that $B=\{\beta_1(t), t\in K\}$ is
reasonably dense in~$\beta_1$, again, for simplicity, we still assume the partition is uniform.
Thus, with $K$ as above, for any $d$, given $t_0\in K$, $R\in SO(d)$, $\gamma\in\Gamma$, we
discretize~\eqref{E:dist-energy} with the trapezoidal rule:
\begin{equation}\label{E:discr-energy}
E^{discr}(t_0,R,\vec{\gamma})=
\sum_{l=1}^{N}h_l'\,\|\sqrt{\dot{\gamma}_l}Rq_1(t_0+\gamma_l)-q_2(t_l)\|^2,
\end{equation}
where $h_1'=(t_2-t_1)/2$, $h_N'=(t_N-t_{N-1})/2$, $h_l'=(t_{l+1}-t_{l-1})/2$ for
$l=2,\ldots,N-1$, $\vec{\gamma}=(\gamma_l)_{l=1}^N$,
$\gamma_1=0, \gamma_N=1$, $\gamma_l = \gamma(t_l)$,
$\dot{\gamma}_l=(\gamma_{l+1} - \gamma_l)/h_l$ for $l=1,\ldots,N-1$,
$\dot{\gamma}_N=\dot{\gamma}_1$,
$h_l = t_{l+1}-t_l$ for $l=1,\ldots,N-1$,
and $q_1(t_0+\gamma_l)$, $l=1,\ldots,N$,
are approximations of $q_1$ at each value of $t_0+\gamma_l$ obtained
from the interpolation of $q_1(t_l)$, $l=1,\ldots,N$, by a cubic~spline.
Accordingly, in \cite{srivastava,srivastava2} \eqref{E:discr-energy} is minimized using a procedure
based on alternating computations of approximately optimal diffeomorphisms (for reparametrizing the
first curve) and approximately optimal rotation matrices (for rotating the first curve) computed
as described in Section 2 and Section~3, respectively, one starting point of $\beta_1$ at
a time. For the sake of completeness, the procedure, with $d$, $K$, $q_1(t_l)$, $q_2(t_l)$,
$l=1,\ldots,N$, as input, is summarized in Procedure~1 below.
We note that in Section~2 and Section~3 computations are carried out that
involve shape functions $q_1$, $q_2$, that if discretized by the same partition $\{t_l\}_{l=1}^N$
are given as finite lists of points: $q_1(t_l)$, $q_2(t_l)$,
$l=1,\ldots,N$.  Thus, in what follows, given shape functions $q$, $\hat{q}$, to say
``\,Execute DP algorithm for $q(t_l)$, $\hat{q}(t_l)$, $l=1,\ldots,N$\," will mean the DP algorithm
(\emph{adapt-DP}) should be executed with $q$, $\hat{q}$ taking the place of $q_1$, $q_2$,
respectively, in
%the computations that are part of
the DP algorithm as outlined in Section~2. The same for the
Kabsch-Umeyama algorithm, i.e., the KU algorithm, as outlined in Section~3, and in the next section
the KU2 algorithm which is the Kabsch-Umeyama algorithm using the~FFT.
%
%\pagebreak
%\begin{algorithm}
%\caption{The optimization algorithm in \cite{srivastava2}}
%\label{A:Srivastava-algo}
\begin{algorithmic}
\STATE \noindent\rule{13cm}{0.4pt}
\STATE {\bf Procedure 1}
\STATE \noindent\rule[.1in]{13cm}{0.4pt}
\STATE Set $h_l=t_{l+1}-t_l$ for $l=1,\ldots,N-1$.
\STATE Set $h_1'=(t_2-t_1)/2$, $h_N'=(t_N-t_{N-1})/2$, $h_l'=(t_{l+1}-t_{l-1})/2$ for
\STATE $l=2,\ldots,N-1$.
\FOR {$\mathrm{each\ }t_0\in K$}
%\STATE {\small\bf if} $t_0=t_1$ {\small\bf then} set $\hat{q}_1(t_l)=q_1(t_l)$ for $l=1,\ldots,N$.
%\STATE {\small\bf else} set $\hat{q}_1(t_l)=q_1(t_0+t_l)$ for $l=1,\ldots,N$. {\small\bf end if}
\STATE Set $\hat{q}_1(t_l)=q_1(t_0+t_l)$ for $l=1,\ldots,N$.
\STATE {\small\bf if} $d\not=1$ {\small\bf then} execute KU algorithm for $q_2(t_l)$,
$\hat{q}_1(t_l)$, $l=1,\ldots,N$,
\STATE to get rotation matrix $R$.
\STATE {\small\bf else} set $R=1$. {\small\bf end if}
\STATE Set $\bar{q}_1(t_l)=R\hat{q}_1(t_l)$ for $l=1,\ldots,N$.
\STATE Execute DP algorithm for $q_2(t_l)$, $\bar{q}_1(t_l)$, $l=1,\ldots,N$, to get
\STATE discretized diffeomorphism $\vec{\gamma} = (\gamma_l)_{l=1}^N$.
\STATE Set $\dot{\gamma}_l=(\gamma_{l+1}-\gamma_l)/h_l$ for $l=1,\ldots,N-1$,
$\dot{\gamma}_N=\dot{\gamma}_1$.
\STATE From interpolation of $q_1(t_l)$, $l=1,\ldots,N$, with a cubic spline\\
       set $\hat{q}_1(t_l)= \sqrt{\dot{\gamma}_l}q_1(t_0+\gamma_l)$ for $l=1,\ldots,N$.
\STATE {\small\bf if} $d\not=1$ {\small\bf then} execute KU algorithm for $q_2(t_l)$,
$\hat{q}_1(t_l)$, $l=1,\ldots,N$,
\STATE to get rotation matrix $R$.
\STATE {\small\bf else} set $R=1$. {\small\bf end if}
\STATE Set $\hat{q}_1(t_l)= R\hat{q}_1(t_l)$ for $l=1,\ldots,N$.
\STATE Compute $E^{discr}(t_0,R,\vec{\gamma})=
\sum_{l=1}^{N}h_l'\,\|\sqrt{\dot{\gamma}_l}Rq_1(t_0+\gamma_l)-q_2(t_l)\|^2$
\STATE $=\sum_{l=1}^{N}h_l'\,\|\hat{q}_1(t_l) - q_2(t_l) \|^2$.
\STATE Keep track of triple $(t_0,R,\vec{\gamma})$, $E^{discr}(t_0,R,\vec{\gamma})$,
\STATE $\hat{q}_1(t_l)$, $l=1,\ldots,N$, with smallest value for $E^{discr}(t_0,R,\vec{\gamma})$.
\ENDFOR
\STATE From interpolation of $\beta_1(t_l)$, $l=1,\ldots,N$, with a cubic spline\\
       set $\hat{\beta}_1(t_l)= R\beta_1(t_0+\gamma_l)$ for $l=1,\ldots,N$.
\STATE Return $(t_0,R,\vec{\gamma})$, $E^{discr}(t_0,R,\vec{\gamma})$, $\hat{\beta}_1(t_l)$,
$\hat{q}_1(t_l)$, $l=1,\ldots,N$.
\STATE \noindent\rule{13cm}{0.4pt}
\end{algorithmic}
%\end{algorithm}
%
On output, $E^{discr}(t_0,R,\vec{\gamma})$ is interpreted to be the elastic shape
distance between $\beta_1$ and $\beta_2$. On the other hand, $\hat{\beta}_1(t_l)$ and $\beta_2(t_l)$,
$l=1,\ldots,N$, are interpreted to achieve the elastic registration of $\beta_1$ and $\beta_2$.
\smallskip\\
Although Procedure~1 above is set up to handle the particular case in which $d$ equals~1 and neither
curve is closed, we note, for the obvious reasons, that if Procedure~1 is adjusted appropriately, the
requirement that the two curves $\beta_1$ and $\beta_2$ must have the same number of points with the
same partition, is then no longer necessary for this case, this case then being the only case in which the
requirement can be ignored. The adjusted Procedure~1, with $q_1(t_l)$,
$l=1,\ldots,N$, $q_2(z_j)$, $j=1,\ldots,M$, as input, is summarized in Procedure~1' below. There
to say ``\,Execute DP algorithm for $q_1(t_l)$, $l=1,\ldots,N$,
$q_2(z_j)$, $j=1,\ldots,M$\," will mean the DP algorithm (\emph{adapt-DP}) should be executed
with $q_1$, $q_2$ exactly as $q_1$, $q_2$ appear in
the DP algorithm as outlined in Section~2.
\begin{algorithmic}
\STATE \noindent\rule{13cm}{0.4pt}
\STATE {\bf Procedure 1'}
\STATE \noindent\rule[.1in]{13cm}{0.4pt}
\STATE Set $h_l=t_{l+1}-t_l$ for $l=1,\ldots,N-1$.
\STATE Set $h_1'=(t_2-t_1)/2$, $h_N'=(t_N-t_{N-1})/2$, $h_l'=(t_{l+1}-t_{l-1})/2$ for
\STATE $l=2,\ldots,N-1$.
\STATE Execute DP algorithm for $q_1(t_l)$, $l=1,\ldots,N$, $q_2(z_j)$, $j=1,\ldots,M$,
\STATE to get discretized diffeomorphism $\vec{\gamma} = (\gamma_l)_{l=1}^N$.
\STATE Set $\dot{\gamma}_l=(\gamma_{l+1}-\gamma_l)/h_l$ for $l=1,\ldots,N-1$,
$\dot{\gamma}_N=\dot{\gamma}_1$.
\STATE From interpolation of $q_2(z_j)$, $j=1,\ldots,M$, with a cubic spline\\
       set $\hat{q}_2(t_l)= \sqrt{\dot{\gamma}_l}q_2(\gamma_l)$ for $l=1,\ldots,N$.
\STATE Compute $E^{discr}(0,1,\vec{\gamma})=
\sum_{l=1}^{N}h_l'\,\|q_1(t_l)-\sqrt{\dot{\gamma}_l}q_2(\gamma_l)\|^2$
\STATE $=\sum_{l=1}^{N}h_l'\,\|q_1(t_l)-\hat{q}_2(t_l)\|^2$.
\STATE From interpolation of $\beta_2(z_j)$, $j=1,\ldots,M$, with a cubic spline\\
       set $\hat{\beta}_2(t_l)= \beta_2(\gamma_l)$ for $l=1,\ldots,N$.
\STATE Return $(0,1,\vec{\gamma})$, $E^{discr}(0,1,\vec{\gamma})$, $\hat{\beta}_2(t_l)$,
$\hat{q}_2(t_l)$, $l=1,\ldots,N$.
\STATE \noindent\rule{13cm}{0.4pt}
\end{algorithmic}
On output, $E^{discr}(0,1,\vec{\gamma})$ is interpreted to be the elastic shape
distance between $\beta_1$ and $\beta_2$. On the other hand, $\beta_1(t_l)$ and $\hat{\beta}_2(t_l)$,
$l=1,\ldots,N$, are interpreted to achieve the elastic registration of $\beta_1$ and $\beta_2$.
\par Finally, with $K$, $q_1(t_l)$, $q_2(t_l)$, $l=1,\ldots,N$, as above, and the partition
$\{t_l\}_{l=1}^N$ uniform if at least one of the curves is closed, given $t_0\in K$, $R\in SO(d)$,
$\gamma\in\Gamma$, we discretize~\eqref{E:dist2-energy} with the trapezoidal rule as follows:
\begin{equation}\label{E:discr2-energy}
E^{discr}(t_0,R,\vec{\gamma})=
\sum_{l=1}^{N}h_l'\,\|Rq_1(t_0+t_l)-\sqrt{\dot{\gamma}_l}Rq_2(\gamma_l)\|^2,
\end{equation}
where $h_1'=(t_2-t_1)/2$, $h_N'=(t_N-t_{N-1})/2$, $h_l'=(t_{l+1}-t_{l-1})/2$ for
$l=2,\ldots,N-1$, $\vec{\gamma}=(\gamma_l)_{l=1}^N$,
$\gamma_1=0, \gamma_N=1$, $\gamma_l = \gamma(t_l)$,
$\dot{\gamma}_l=(\gamma_{l+1} - \gamma_l)/h_l$ for $l=1,\ldots,N-1$,
$\dot{\gamma}_N=\dot{\gamma}_1$,
$h_l = t_{l+1}-t_l$ for $l=1,\ldots,N-1$,
and $q_2(\gamma_l)$, $l=1,\ldots,N$,
are approximations of $q_2$ at each $\gamma_l$ obtained from the interpolation of $q_2(t_l)$,
$l=1,\ldots,N$, by a cubic~spline.  Accordingly, for the purpose of minimizing \eqref{E:discr2-energy},
we use a procedure that alternates computations, as described in Section~2 and Section~3, of
approximately optimal diffeomorphisms (a constant number of them per iteration for reparametrizing the
second curve) and successive computations of approximately optimal rotation matrices (for rotating the
first curve) for all starting points of the first curve. As noted in~\cite{dogan2}, carrying out
computations this way is not only more efficient all by itself, but, if both curves are closed, allows
applications of the Fast Fourier Transform (FFT) as demonstrated in \cite{dogan} for~$d=2$, for
computing successively in an even more efficient manner, as described in the next section, optimal
rotation matrices for all starting points of the first~curve. The procedure, with $K$, $q_1(t_l)$,
$q_2(t_l)$, $l=1,\ldots,N$, as input, is summarized in Procedure~2 below. Note, $itop$ in the
procedure is an input variable that must be set equal to a positive integer that is constant relative
to~$N$, and not larger than the cardinality of~$K$. It is the number of times the second {\small\bf for}
loop in the {\small\bf repeat} loop of the procedure is executed during each iteration of the {\small\bf repeat}
loop. It is in the second {\small\bf for} loop that the DP algorithm is executed, thus the execution time
of the procedure can be large if $itop$ is greater than~1. Actually the first {\small\bf for} loop in the
{\small\bf repeat} loop of the procedure takes a lot less time than the second {\small\bf for} loop
even if $itop$ equals~1. It is in the first {\small\bf for} loop that the KU algorithm is executed. We note,
in our experiments, $itop$ equal to~1 has usually sufficed for curves of relatively simple curvature, e.g.,
spherical ellipsoids in $3-$dimensional space (see Section~6). For curves of more complex curvatures, higher
values have usually been required for the successful execution of the procedure.
%For more complex curves,
%\pagebreak
%
\begin{algorithmic}
\STATE \noindent\rule{13cm}{0.4pt}
\STATE {\bf Procedure 2}
\STATE \noindent\rule[.1in]{13cm}{0.4pt}
\STATE Set $h_l=t_{l+1}-t_l$ for $l=1,\ldots,N-1$.
\STATE Set $h_1'=(t_2-t_1)/2$, $h_N'=(t_N-t_{N-1})/2$, $h_l'=(t_{l+1}-t_{l-1})/2$ for
\STATE $l=2,\ldots,N-1$.
\STATE Set $\hat{q}_2(t_l)=q_2(t_l)$ for $l=1,\ldots,N$.
\STATE Set $iter=0$, $E^{curr}=10\,^6$, $iten=10$, $tol=10\,^{-6}$.
\REPEAT
\STATE Set $iter=iter+1$, $E^{prev} = E^{curr}$.
\FOR {$\mathrm{each\ }t_0\in K$}
%\STATE {\small\bf if} $t_0=t_1$ {\small\bf then} set $\hat{q}_1(t_l)=q_1(t_l)$ for $l=1,\ldots,N$.
%\STATE {\small\bf else} set $\hat{q}_1(t_l)=q_1(t_0+t_l)$ for $l=1,\ldots,N$. {\small\bf end if}
\STATE Set $\hat{q}_1(t_l)=q_1(t_0+t_l)$ for $l=1,\ldots,N$.
\STATE Execute KU algorithm for $\hat{q}_2(t_l)$, $\hat{q}_1(t_l)$, $l=1,\ldots,N$, to get
rotation matrix $R$ and $maxtrace$.
\STATE Identify $(t_0,R)$ as a couple of interest and associate with it the value of $maxtrace$.
\STATE Keep track of identified couples of interest $(t_{0i},R_i)$, $i=1,\ldots,itop$, satisfying that
for each $i$, $i=1,\ldots,itop$, the value of $maxtrace$ associated with $(t_{0i},R_i)$ is one of the
$itop$ largest values among the values of $maxtrace$ associated with all couples of interest identified
so far.
\ENDFOR
%\STATE {\small\bf if} $t_0=t_1$ {\small\bf then} set $\hat{q}_1(t_l)=Rq_1(t_l)$ for $l=1,\ldots,N$.
%\STATE {\small\bf else} set $\hat{q}_1(t_l)=Rq_1(t_0+t_l)$ for $l=1,\ldots,N$. {\small\bf end if}
\FOR {$i=1,\ldots,itop$}
\STATE Set $t_0=t_{0i}$, $R=R_i$.
\STATE Set $\hat{q}_1(t_l)=Rq_1(t_0+t_l)$ for $l=1,\ldots,N$.
\STATE Execute DP algorithm for $\hat{q}_1(t_l)$, $q_2(t_l)$, $l=1,\ldots,N$, to get
\STATE discretized diffeomorphism $\vec{\gamma} = (\gamma_l)_{l=1}^N$.
\STATE Set $\dot{\gamma}_l=(\gamma_{l+1}-\gamma_l)/h_l$ for $l=1,\ldots,N-1$,
       $\dot{\gamma}_N=\dot{\gamma}_1$.
\STATE From interpolation of $q_2(t_l)$, $l=1,\ldots,N$, with a cubic spline\\
       set $\hat{q}_2(t_l)= \sqrt{\dot{\gamma}_l}q_2(\gamma_l)$ for $l=1,\ldots,N$.
\STATE Compute $E^{curr}=E^{discr}(t_0,R,\vec{\gamma})=
\sum_{l=1}^{N}h_l'\,\|Rq_1(t_0+t_l)-\sqrt{\dot{\gamma}_l}q_2(\gamma_l)\|^2$
\STATE $=\sum_{l=1}^{N}h_l'\,\|\hat{q}_1(t_l) - \hat{q}_2(t_l) \|^2$.
\STATE Keep track of triple $(t_0,R,\vec{\gamma})$, $E^{curr}$,
$\hat{q}_1(t_l)$, $\hat{q}_2(t_l)$, $l=1,\ldots,N$, with smallest value for~$E^{curr}$.
\ENDFOR
\UNTIL {$|E^{curr}-E^{prev}|< tol$ or $iter > iten$.}
\STATE From interpolation of $\beta_2(t_l)$, $l=1,\ldots,N$, with a cubic spline\\
       set $\hat{\beta}_2(t_l)= \beta_2(\gamma_l)$ for $l=1,\ldots,N$.
%\STATE {\small\bf if} $t_0=t_1$ {\small\bf then} set $\hat{\beta}_1(t_l)=R\beta_1(t_l)$ for $l=1,\ldots,N$.
%\STATE {\small\bf else} set $\hat{\beta}_1(t_l)=R\beta_1(t_0+t_l)$ for $l=1,\ldots,N$. {\small\bf end if}
\STATE Set $\hat{\beta}_1(t_l)=R\beta_1(t_0+t_l)$ for $l=1,\ldots,N$.
\STATE Return $(t_0,R,\vec{\gamma})$, $E^{discr}(t_0,R,\vec{\gamma})\ (=E^{curr})$,
$\hat{\beta}_1(t_l)$, $\hat{q}_1(t_l)$, $\hat{\beta}_2(t_l)$, $\hat{q}_2(t_l)$,\\ $l=1,\ldots,N$.
\STATE \noindent\rule{13.2cm}{0.4pt}
\end{algorithmic}
On output, $E^{discr}(t_0,R,\vec{\gamma})$ is interpreted to be the elastic shape
distance between $\beta_1$ and $\beta_2$. On the other hand, $\hat{\beta}_1(t_l)$ and
$\hat{\beta}_2(t_l)$, $l=1,\ldots,N$, are interpreted to achieve the elastic registration of
$\beta_1$ and $\beta_2$.
\smallskip\\
Similar to Procedure~1, Procedure~2 above is also set up to handle the particular case in
which $d$ equals~1 and neither curve is closed. But since similar to Procedure~1 the requirement
that the two curves $\beta_1$ and $\beta_2$ must have the same number of points with the same
partition, is not necessary for this case, Procedure~1' can be used instead of Procedure~$2$
in the absence of the requirement.
%\section{\large Successive Computations of Optimal Rotations
%in $d-$dimensional Space with FFT}
\section{\large Successive Computations of Rotations with FFT for Rigid alignment
of Curves in $d-$dimensional Space}
Again, let $\beta_1$, $\beta_2$, $q_1$, $q_2$ be as above, i.e.,
$\beta_n: [0,1]\rightarrow \mathbb{R}^d$, $n=1,2$, are absolutely continuous functions representing
simple curves in~$\mathbb{R}^d$ of unit length, and $q_n: [0,1]\rightarrow \mathbb{R}^d$, $n=1,2$, are
square-integrable functions that are the shape functions or SRVF's of $\beta_n$, $n=1,2$, respectively.
In this section, using arguments similar to those used in \cite{dogan} for $d=2$, we first present an
alternative version of the KU algorithm that uses the FFT for the
purpose of speeding up the successive computations in Procedure~$2$ in the previous section, of
approximately optimal rotation matrices for all starting points of one of the curves. These
computations actually take place in the first {\bf for} loop of that procedure. Taking into account the
nature of the FFT, we assume both curves are closed (this will become evident below), and without any
loss of generality, for the purpose of developing the alternative version of the KU algorithm in a
manner similar to the way in which the KU algorithm was developed in Section~3, assume that any point
in $\beta_2$ can be treated as a starting point of~$\beta_2$, and that it is $\beta_1$ that has a fixed
starting point. In particular, it follows then that $q_2$ can be interpreted to be a periodic function
from $\mathbb{R}$ into $\mathbb{R}^d$, $q_2(t+1)=q_2(t)$ for all vaues of~$t$. Taking into account as
well the part of Procedure~$2$ in the previous section that we are trying to improve (the first {\bf for}
loop), and following the reasoning in Section~3 to obtain~\eqref{E:energy0}, ideally,
we would like to solve a problem of the following type: Find $t_0\in [0,1]$, and a $d\times d$
rotation matrix~$R$ that minimize
\begin{equation}\label{E:energy01}
E(t_0,R) = \int_0^1 \| q_1(t) - R q_2(t_0+t) \|^2 dt.
\end{equation}
As noted above, in practice, we work with curves $\beta_1$, $\beta_2$, given as discrete sets of
points. Accordingly, for some integer $N>0$, and a partition of~$[0,1]$, $\{t_l\}_{l=1}^{N}$,
$t_1=0<t_2<\ldots<t_N=1$, for $n=1,2$, the curve $\beta_n$ is given as a list
of $N$ points in the curve, where for $l=1,\ldots,N$, $\beta_n(t_l)$ is the $l^{th}$ point in the
list for~$\beta_n$. Similarly for $q_1$, $q_2$, except that for $l=1,\ldots,N$, $q_1(t_l)$ and
$q_2(t_l)$ are computed as described in Section~2. Again, we assume $K$ as defined in the previous
section equals $\{t_1,\ldots,t_{N-1}\}$ or a subset of it, a subset essentially of size~$O(N)$,
so that $\{\beta_2(t), t\in K\}$ is then interpreted to be the set of starting points
of~$\beta_2$. Also, as justified in the previous section, the partition $\{t_l\}_{l=1}^N$ must then
be uniform. In what follows, for $l=1,\ldots,N$, $k=1,\ldots,d$, $j=1,\ldots,d$, $q_1^l$ is
$q_1(t_l)$, $q_2^l$ is $q_2(t_l)$, $q_{1k}^l$ is the $k^{th}$ coordinate of $q_1^l$, $q_{2j}^l$ is
the $j^{th}$ coordinate of $q_2^l$, and $\hat{q}\,^l_{1k}$ is $q\,^{N-\,l+1}_{1k}$.
\par
In order to discretize integral \eqref{E:energy01}, we define for
each $m=1,\ldots,N-1$, points $q_2^{m\,\oplus\,l}, l=1,\ldots,N$, by
\begin{equation*}
q_2^{m\,\oplus\,l} = q_2(t_m+t_l),
\end{equation*}
and let $q_{21}^{m\,\oplus\,l},\ldots,q_{2d}^{m\,\oplus\,l}$ be the
$d$ coordinates of $q_2^{m\,\oplus\,l}$ so that
\begin{equation*}
(q_{21}^{m\,\oplus\,l},\ldots,q_{2d}^{m\,\oplus\,l})^T=q_2^{m\,\oplus\,l}.
\end{equation*}
We note as well that for $m=1,\ldots,N-1$, we may then assume the
existence of additional functions $q_2^m: [0,1]\rightarrow \mathbb{R}^d$, given in their
discretized form as
\begin{equation*}
q_2^m(t_l)=q_2(t_m+t_l)=q_2^{m\,\oplus\,l}, l=1,\ldots,N.
\end{equation*}
With $1\leq m\leq N-1$, letting $h=1/(N-1)$, we then discretize integral \eqref{E:energy01}
using the uniform trapezoidal rule for when both curves are closed:
\begin{equation}\label{E:energy02}
E^{discr}(m,R) = h \sum_{l=1}^{N-1} \| q_1(t_l) - R q_2^m(t_l) \|^2
= h \sum_{l=1}^{N-1} \| q_1^l - R q_2^{m\,\oplus\,l} \|^2.
\end{equation}
%Given an integer $m$, $1\leq m\leq N-1$, and a $d\times d$ rotation matrix $R$ that
%minimizes~\eqref{E:energy02}, we identify $(m,R)$ as a couple of interest
%and associate with it the corresponding value of~\eqref{E:energy02}.
%Thus, the problem of finding $t_0\in [0,1]$ and a $d\times d$ rotation matrix
%$R$ that minimize~\eqref{E:energy01} becomes the problem of finding for each~$m$, $m=1,\ldots,N-1$,
%a $d\times d$ rotation matrix $R$ that minimizes~\eqref{E:energy02}, identifying $(m,R)$ as a
%couple of interest, associating with it the corresponding value of~\eqref{E:energy02}, and
%finally identifying an integer $m$, \mbox{$1\leq m\leq N-1$}, and a $d\times d$ rotation matrix $R$
%such that $(m,R)$ is one of the identified couples of interest and the value of~\eqref{E:energy02}
%associated with it is the smallest among the values of~\eqref{E:energy02} associated with any of
%the identified couples of interest.\\
Thus, the problem of finding $t_0\in [0,1]$ and a $d\times d$ rotation matrix $R$ that
minimize~\eqref{E:energy01} becomes the problem of finding~$m$, $1\leq m\leq N-1$, with~$t_m$
in~$K$, and a $d\times d$ rotation matrix $R$ that minimize~\eqref{E:energy02}.\\
For this purpose, for each $m$, $m=1,\ldots,N-1$, we define a $d\times d$ matrix $A(m)$ by defining
its entries $A_{kj}(m)$ for each pair $k,j = 1,\ldots, d$, by
\begin{equation}\label{A:entries}
A_{kj}(m) = \sum_{l=1}^{N-1} q_{1k}^l q_{2j}^{m\,\oplus\,l},
\end{equation}
so that for fixed $m$, minimizing \eqref{E:energy02} over all $d\times d$ rotation matrices $R$ is
equivalent to maximizing
\begin{equation}\label{E:disc-energy02}
\sum_{l=1}^{N-1} (q_1^l)^T R q_2^{m\,\oplus\,l} = \mathrm{tr}(R A(m)^T).
\end{equation}
We note that for fixed $m$, we can execute the KU algorithm for $q_1(t_l)$, $q_2^m(t_l)$,
$l=1,\ldots,N$, to compute $R$ that maximizes~\eqref{E:disc-energy02}.
Doing this for each $m$, $1\leq m\leq N-1$, with~$t_m$ in~$K$, we identify among them an $m$ for
which the maximization of~\eqref{E:disc-energy02} is the largest. The solution is then that $m$
together with the rotation matrix $R$ at which the maximization is achieved. We also note that
computing $A(m)$ for each $m$ is $O(N)$ so that computing $O(N)$ of them is then $O(N^2)$ if each
$A(m)$ is computed separately. This is exactly how it is done in Procedure~$2$ in the previous section.
\par
For each pair $k,j=1,\ldots,d$, with $A_{kj}(m)$ as in~\eqref{A:entries}, we propose to compute
all of $A_{kj}(1),\ldots,A_{kj}(N-1)$ in $O(N\log N)$ time using the FFT to accomplish the Discrete
Fourier Transform (DFT). For this purpose, for $k=1,\ldots,d$, let
$\hat{q}_{1k}=(\hat{q}\,^1_{1k},\ldots,\hat{q}\,^{N-1}_{1k})$, and for $j=1,\ldots,d$, let
$q_{2j}=(q\,^1_{2j},\ldots,q\,^{N-1}_{2j})$.
Given arbitrary vectors $x$, $y$ of length~$N-1$, we let {\bf DFT}$(x)$ and
{\bf DFT}$^{-1}(y)$ denote the DFT of~$x$ and the inverse DFT of~$y$, respectively.
With the symbol $\cdot$ indicating component by component multiplication of two vectors,
then by the convolution theorem for the DFT we have for each pair $k,j = 1,\ldots,d$,
\begin{eqnarray*}
(A_{kj}(1),\ldots,A_{kj}(N-1)) &=&
(\sum_{l=1}^{N-1} q_{1k}^lq_{2j}^{1\,\oplus\,l},\ldots,
 \sum_{l=1}^{N-1} q_{1k}^lq_{2j}^{(N-1)\,\oplus\,l})\\
&=&\mathrm{\bf DFT}^{-1}[\mathrm{\bf DFT}(\hat{q}_{1k})\cdot
\mathrm{\bf DFT}(q_{2j})]
\end{eqnarray*}
which for each pair $k,j=1,\ldots,d$, enables us to reduce the computation of
all of $A_{kj}(1),\ldots,A_{kj}(N-1)$ to three $O(N\log N)$ FFT operations.
Thus, we can compute all of $A(1),\ldots,A(N-1)$ in $O(N\log N)$ time with the~FFT.
\par
An outline of the alternative version of the KU algorithm, the KU2 algorithm, that uses
the FFT, follows.  Here for arbitrary vectors $x$, $y$ of length~$N-1$, {\bf FFT}$(x)$,
{\bf IFFT}$(y)$ denote {\bf DFT}$(x)$, {\bf DFT}$^{-1}(y)$, respectively, computed with
the~FFT. Note, $itop$ in the algorithm is an input variable as described in the previous
section before the outline of Procedure~$2$, an input variable used there exclusively in
that procedure, its purpose to control the number of times the DP algorithm is executed
in the procedure. Here, with the same purpose, before it is an input variable of the KU2
algorithm, it is first an input variable of a
procedure in which the KU2 and DP algorithms are executed, \mbox{Procedure~3},
the outline of which appears later in this section, its purpose discussed as~well.
\begin{algorithmic}
\STATE \noindent\rule{13cm}{0.4pt}
\STATE {\bf Algorithm Kabsch-Umeyama} with FFT (KU2 algorithm)
\STATE \noindent\rule[.1in]{13cm}{0.4pt}
\STATE Set $q^l_{1k}$ equal to $k^{th}$ coordinate of $q_1(t_l)$
for $l=1,\ldots,N$, $k=1,\ldots,d$.
\STATE Set $q^l_{2j}$ equal to $j^{th}$ coordinate of $q_2(t_l)$
for $l=1,\ldots,N$, $j=1,\ldots,d$.
\STATE Set $\hat{q}\,^l_{1k}=q\,^{N-\,l+1}_{1k}$ for $l=1,\ldots,N$, $k=1,\ldots,d$.
\STATE Set $\hat{q}_{1k}=(\hat{q}\,^1_{1k},\ldots,\hat{q}\,^{N-1}_{1k})$ for $k=1,\ldots,d$.
\STATE Set $q_{2j}=(q\,^1_{2j},\ldots,q\,^{N-1}_{2j})$, for $j=1,\ldots,d$.
\FOR {each pair $k,j=1,\ldots,d$}
\STATE Compute $(A_{kj}(1),\ldots,A_{kj}(N-1))=$ {\bf IFFT}[{\bf FFT} $(\hat{q}_{1k})\,\cdot\,$
              {\bf FFT}$(q_{2j})$].
\ENDFOR
\FOR {$\mathrm{each\ }m$, $1\leq m \leq N-1$, with $t_m\in K$}
\STATE Identify $d\times d$ matrix $A(m)$ with entries $A_{kj}(m)$ for\\ each pair $k,j=1,\ldots,d$.
\STATE Compute SVD of $A(m)$, i.e., identify $d\times d$ matrices $U$, $S$, $V$,
\STATE so that $A(m) = U S V^T$ in the SVD sense.
\STATE Set $s_1= \ldots = s_{d-1}=1$.
\STATE {\small\bf if} $\det(UV) > 0$ {\small\bf then} set $s_d=1$.
\STATE {\small\bf else} set $s_d=-1$. {\small\bf end if}
\STATE Set $\tilde{S} = \mathrm{diag}\{s_1,\ldots,s_d\}$.
\STATE Compute $R = U \tilde{S} V^T$ and $maxtrace = \mathrm{tr}(RA(m)^T)$.
\STATE Identify $(m,R)$ as a couple of interest and associate with it the value of $maxtrace$.
\STATE Keep track of identified couples of interest $(m_i,R_i)$, $i=1,\ldots,itop$, satisfying that
for each $i$, $i=1,\ldots,itop$, the value of $maxtrace$ associated with $(m_i,R_i)$ is one of the
$itop$ largest values among the values of $maxtrace$ associated with all couples of interest identified
so far.
\ENDFOR
\STATE Return couples $(m_i,R_i)$, $i=1,\ldots,itop$.
\STATE \noindent\rule{13cm}{0.4pt}
\end{algorithmic}
\smallskip
\par We note that if $d=1$, the KU2 algorithm still computes couples $(m_i,R_i)$, $i=1,\ldots,itop$,
with the resulting $R_i$'s always equal to~1.
\\ \smallskip
\par Finally, a modified version of Procedure~$2$ in the previous section, Procedure~3, follows.
Here $\beta_1$, $\beta_2$, $q_1$, $q_2$, $N$, $\{t_l\}_{l=1}^N$, $\beta_1(t_l)$, $\beta_2(t_l)$,
$q_1(t_l)$, $q_2(t_l)$, $l=1,\ldots,N$, are as above. Thus, $\beta_1$, $\beta_2$ are closed
curves and $\{t_l\}_{l=1}^N$ is uniform.
The procedure with $K$, $itop$, $q_1(t_l)$, $q_2(t_l)$, $l=1,\ldots,N$, as input, is essentially the
same as Procedure~$2$, except that the {\bf for} loop in Procedure~2 that executes the KU~algorithm
of Section~3 as many times as there are starting points of the first curve ($\beta_1$), is replaced
by the execution of the KU2~algorithm outlined above. This has the effect of speeding up the successive
computations appearing in Procedure~$2$ of approximately optimal rotation matrices for all starting
points of the first curve due to the fact that the KU2 algorithm uses the FFT which takes $O(N\log N)$
time, while the {\bf for} loop in Procedure~2 computes each approximately optimal rotation matrix
separately thus taking $O(N^2)$~time.
%thus using the FFT for the purpose of speeding
%up the successive computations in Procedure~$2$ of approximately optimal rotation matrices
%for all starting points of the first curve~($\beta_1$).
%
\begin{algorithmic}
\STATE \noindent\rule{13cm}{0.4pt}
\STATE {\bf Procedure 3}
\STATE \noindent\rule[.1in]{13cm}{0.4pt}
\STATE Set $h=1/(N-1)$.
\STATE Set $\hat{q}_2(t_l)=q_2(t_l)$ for $l=1,\ldots,N$.
\STATE Set $iter=0$, $E^{curr}=10\,^6$, $iten=10$, $tol=10\,^{-6}$.
\REPEAT
\STATE Set $iter=iter+1$, $E^{prev} = E^{curr}$.
\STATE Execute KU2 algorithm for $\hat{q}_2(t_l)$, $q_1(t_l)$, $l=1,\ldots,N$, to get
\STATE couples $(m_i,R_i)$, $m_i$ an integer, $1\leq m_i\leq N-1$, $R_i$ a rotation matrix,
$i=1,\ldots,itop$.
\FOR {$i=1,\ldots,itop$}
\STATE Set $m=m_i$, $R=R_i$.
\STATE Set $\hat{q}_1(t_l)=Rq_1(t_m+t_l)$ for $l=1,\ldots,N$.
\STATE Execute DP algorithm for $\hat{q}_1(t_l)$, $q_2(t_l)$, $l=1,\ldots,N$, to get
\STATE discretized diffeomorphism $\vec{\gamma} = (\gamma_l)_{l=1}^N$.
\STATE Set $\dot{\gamma}_l=(\gamma_{l+1}-\gamma_l)/h$ for $l=1,\ldots,N-1$,
       $\dot{\gamma}_N=\dot{\gamma}_1$.
\STATE From interpolation of $q_2(t_l)$, $l=1,\ldots,N$, with a cubic spline\\
       set $\hat{q}_2(t_l)= \sqrt{\dot{\gamma}_l}q_2(\gamma_l)$ for $l=1,\ldots,N$.
\STATE Compute $E^{curr}=E^{discr}(t_m,R,\vec{\gamma})$
\STATE $=\sum_{l=1}^{N-1}h\,\|Rq_1(t_m+t_l)-\sqrt{\dot{\gamma}_l}q_2(\gamma_l)\|^2$
$=\sum_{l=1}^{N-1}h\,\|\hat{q}_1(t_l) - \hat{q}_2(t_l) \|^2$.
\STATE Keep track of triple $(t_m,R,\vec{\gamma})$, $E^{curr}$,
$\hat{q}_1(t_l)$, $\hat{q}_2(t_l)$, $l=1,\ldots,N$, with smallest value for~$E^{curr}$.
\ENDFOR
\UNTIL {$|E^{curr}-E^{prev}|< tol$ or $iter > iten$.}
\STATE From interpolation of $\beta_2(t_l)$, $l=1,\ldots,N$, with a cubic spline\\
       set $\hat{\beta}_2(t_l)= \beta_2(\gamma_l)$ for $l=1,\ldots,N$.
%\STATE {\small\bf if} $t_m=t_1$ {\small\bf then} set $\hat{\beta}_1(t_l)=R\beta_1(t_l)$ for $l=1,\ldots,N$.
%\STATE {\small\bf else} set $\hat{\beta}_1(t_l)=R\beta_1(t_m+t_l)$ for $l=1,\ldots,N$. {\small\bf end if}
\STATE Set $\hat{\beta}_1(t_l)=R\beta_1(t_m+t_l)$ for $l=1,\ldots,N$.
\STATE Return $(t_m,R,\vec{\gamma})$, $E^{discr}(t_m,R,\vec{\gamma})\ (=E^{curr})$,
$\hat{\beta}_1(t_l)$, $\hat{q}_1(t_l)$, $\hat{\beta}_2(t_l)$, $\hat{q}_2(t_l)$,\\ $l=1,\ldots,N$.
\STATE \noindent\rule{13.2cm}{0.4pt}
\end{algorithmic}
On output, $E^{discr}(t_m,R,\vec{\gamma})$ is interpreted to be the elastic shape
distance between $\beta_1$ and $\beta_2$. On the other hand, $\hat{\beta}_1(t_l)$ and
$\hat{\beta}_2(t_l)$, $l=1,\ldots,N$, are interpreted to achieve the elastic registration of
$\beta_1$ and $\beta_2$.
\section{\large Results from Computations with Implementation of Methods}
A software package that incorporates the methods presented in this paper for computing the
elastic registration of two simple curves in $d-$dimensional space, $d$ a positive integer, and
therefore the elastic shape distance between them, has been implemented. The implementation is in
Matlab\footnote{The identification of any commercial product or trade name does not imply
endorsement or recommendation by the National Institute of Standards and Technology.}
with the exception of the Dynamic Programming routine which is written in Fortran but is executed
as a Matlab mex file. In this section, we present results obtained from executions of the software
package with~$d=3$. We note, the sofware package as well as input data files, a README file, etc.
can be obtained at the following link\smallskip\\
\hspace*{.35in}\verb+https://doi.org/10.18434/mds2-2329+
%\\
%\hspace*{.35in}\verb+http://math.nist.gov/~JBernal+ \verb+/ESD_alldim.zip+
\smallskip\\
With the exception of the Matlab driver routine, ESD\_\,driv\_\,3\,dim.m, which is designed for the
case~$d=3$, all Matlab routines in the package can be executed for any~$d$ if the current driver routine
is adjusted or replaced to handle the value of~$d$. However, parameter dimx in the Fortran routine
DP\_\,MEX\,\_\,WNDSTRP\_\,ALLDIM.F may have to be modified so that instead of having a value of 3, it has
the value of~$d$.  The Fortran routine must then be processed to obtain a new mex file for the routine by
typing in the Matlab window:
mex -\,compatibleArrayDims DP\_\,MEX\,\_\,WNDSTRP\_\,ALLDIM.F
%\smallskip\\
\par
Given discretizations of two simple curves $\beta_1$, $\beta_2$, $\beta_1: [0,T_1]\rightarrow \mathbb{R}^d$,
$\beta_2: [0,T_2]\rightarrow \mathbb{R}^d$, $T_1,\;T_2>0$, the elastic registration of $\beta_1$ and $\beta_2$
to be computed together with the elastic shape distance between them, irrespective of the value of $d$, the
program always proceeds first to compute an approximation of the length of each curve by computing the length
of each line segment joining consecutive points on the curve in the discretization of the curve and adding
these lengths, and then proceeds to scale the two curves so that each curve has approximate length~1 (each
point in the discretization of each curve is divided by the approximate length of the curve). The program
then proceeds to scale, if any, the two partitions that discretize the curves so that they become partitions
of~$[0,1]$, or if no partitions are given, to create two partitions of~$[0,1]$, one for each curve, according to
the number of points in the discretization of each curve, the discretization of each curve then assumed to be the
result of discretizing the curve by the corresponding partition, each partition uniform if at least one curve is
closed, each partition parametrizing the corresponding curve by arc length otherwise. Utilizing the given or
created partitions and the discretizations of the curves, with the exception of the case in which $d$ equals~1
and both curves are open, the program then proceeds to create a common partition of~$[0,1]$ for the two curves
and to discretize each curve by this common partition using cubic splines. If at least one curve is closed, and the
numbers of points in the first curve and second curve are $N$ and $M$, respectively, letting $L$ equal the larger
of $N$ and $M$, the common partition is then taken to be the uniform partition of~$[0,1]$ of size equal to~$L$.
This is in accordance with the requirement established in Section~4 that if at least one curve is closed (the set
of starting points of one of the curves will have more than one point), in order to compute the elastic
shape distance and registration in the appropriate manner, the curves should be discretized by the same uniform
partition. Note that a set of starting points of one of the curves having more than one point is then identified
satisfying that it is the discretization by the uniform partition of one of the curves (a closed curve), or a subset
of it. On the other hand, if both curves are open, $d$ not equal to~1,
the common partition is then taken to be the union of the two partitions discretizing
the curves minus certain points in this union that are eliminated systematically so that the distance between
any two consecutive points in the common partition does not exceed some tolerance. This is in accordance with
the requirement established in Section~4 that if both curves are open (each curve has exactly one starting point),
$d$ not equal to~1, in order to compute the elastic shape distance and registration
in the appropriate manner, the curves should be discretized by the same partition, a partition not
necessarily uniform. Finally, if both curves are open and $d=1$, no common partition is created and the curves
continue to be discretized by the same given or created partitions. All of the above is accomplished by
Matlab routine ESD\_\,comp\_\,alldim.m during the execution of the software package. Once this routine
is done, the actual computations of the elastic shape distance and registration are carried out by Matlab
routine ESD\_\,core\_\,alldim.m in which all of the procedures presented in Section~4 and Section~5 have
been implemented.
\smallskip
\par The results that follow were obtained from applications of our software package on
discretizations of curves in $3-$dimensional space of the helix and spherical ellipsoid kind. With an
observer at the origin of the $3-$dimensional Euclidean space whose line of sight is the positive
$z-$axis, given $T>0$, a circular helix of radius~1 with axis of rotation the positive $z-$axis and
that moves away from the observer in a clockwise screwing motion, is defined by
\[ x(t)=\cos t,\ \ \ y(t) = \sin t,\ \ \ z(t) = t,\ \ \ t\in [0,T]. \]
On the other hand, with an observer at the origin of the $3-$dimensional Euclidean space whose line of sight
is the positive $z-$axis, given $r, a, b$, with $r > a > 0$, $r > b > 0$, a spherical ellipsoid with axis of
rotation the positive $z-$axis and that as viewed by the observer
moves around its axis of rotation in a clockwise direction, is defined by
\[ x(t)=a\cos t,\ \ \ y(t) = b\sin t,\ \ \ z(t) = (r^2-x(t)^2-y(t)^2)^{1/2},\ \ \ t\in [0,2\pi]. \]
We note that in the obvious similar manner, helices and spherical ellipsoids with axis of rotation the
positve/negative $x-,y-,z-$axis can be defined as~well. We also note that as defined above helices are
open curves, and spherical ellipsoids are closed curves.
\par
\begin{figure}
\centering
\begin{tabular}{ccc}
\includegraphics[width=0.3\textwidth]{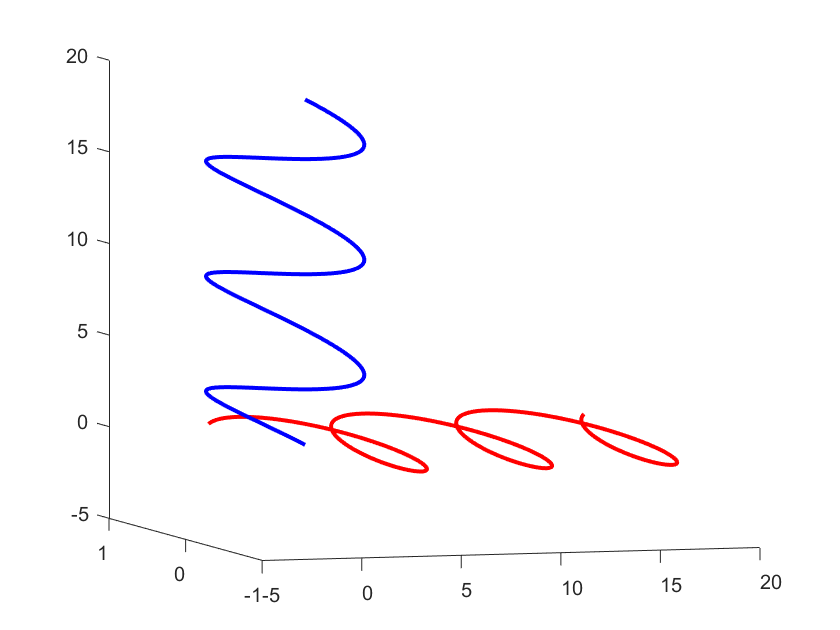}
&
\includegraphics[width=0.3\textwidth]{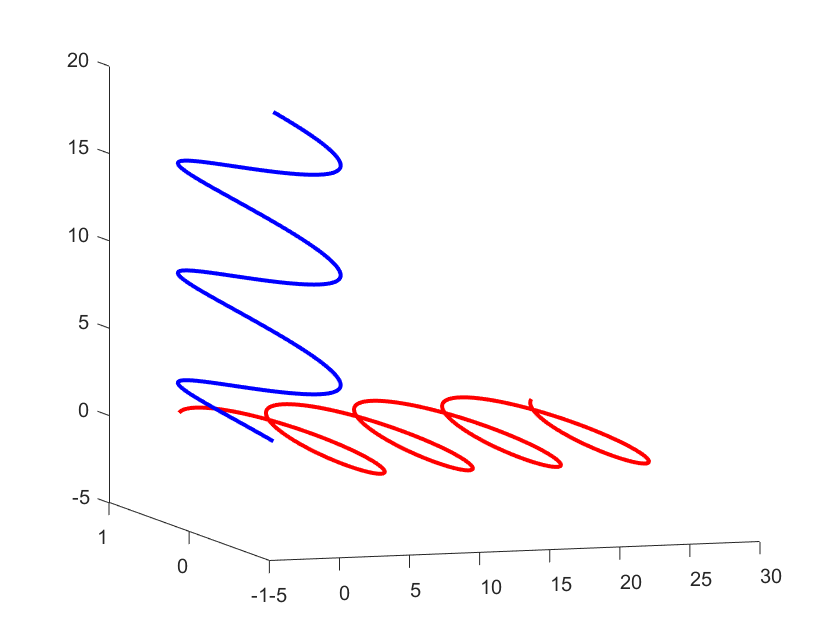}
&
\includegraphics[width=0.3\textwidth]{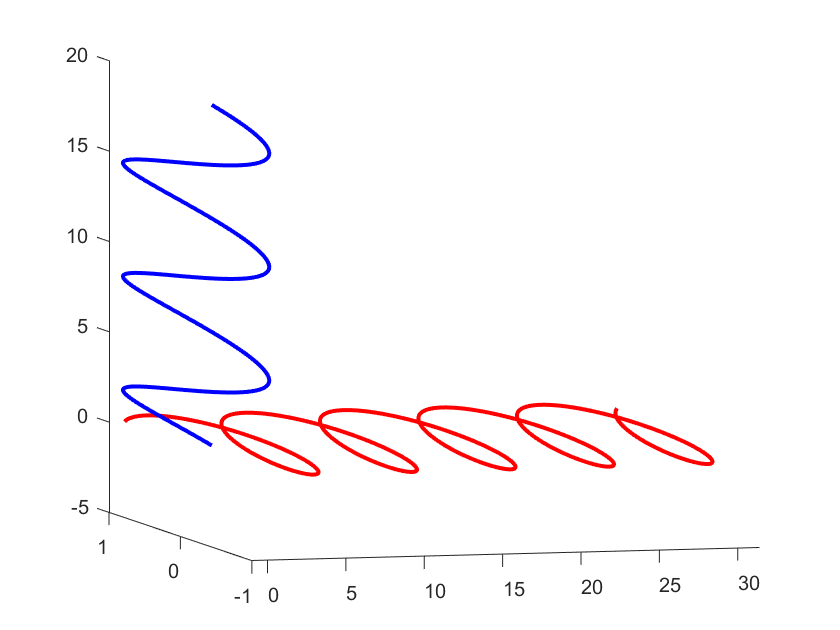}
\end{tabular}
\caption{\label{F:curves}
Three plots of helices. The elastic registration of the two helices
in each plot and the elastic shape distance between them were computed.
}

\end{figure}

 Three plots depicting helices are shown in Figure~\ref{F:curves}. (Note that in the plots there,
the $x-$axis, the $y-$axis and the $z-$axis are not always to scale relative to one another).
In each plot two helices appear. The helix in each plot with the positive $z-$axis as its axis of
rotation was considered to be the first curve or helix in the plot. In each plot this helix was obtained
by setting $T$ to $6\pi$ in the definiton of a helix above so that it has three loops in each plot and thus
is the same helix in all three plots. The other helix in each plot has the positive $x-$axis as its axis
of rotation and was considered to be the second curve or helix in each plot. From left to right in the three
plots, the second helix was obtained by setting $T$ to $6\pi$, $8\pi$, $10\pi$, respectively, in the definition
of a helix above, the definition modified in the obvious manner so that the helix has the positive $x-$axis
as its axis of rotation. Thus the second helix has three, four, five loops, from left to right in the three
plots. All helices in the plots were then discretized as described below and the elastic
registration of the two helices in each plot and the elastic shape distance between them were then computed
through executions of our software package (mostly executions of Procedure~$2$ in Section~4).
Accordingly, one would expect the elastic shape distances, if given in the order of the plots from left to
right, to have been in strictly increasing order with the first distance essentially equal to zero. That is
exactly what we obtained:
\ \ \ 0.00000\ \ \ \ 0.48221\ \ \ \ 0.60352.
\smallskip\\
We note that on input the first helix in each plot was given as the same discretization of a $3-$loop helix,
a helix discretized by a uniform partition of~$[0,6\pi]$, the discretization consisting of 451 points.
On the other hand, the second helix in each plot was given as well as the discretization of a helix, from left
to right in the three plots a helix having 3, 4,~5 loops, respectively, a helix discretized by a uniform partition
of~$[0,6\pi]$, $[0,8\pi]$, $[0,10\pi]$, respectively, the discretizations consisting of 451, 601,~751 points,
respectively. Given a pair of helices in one of the three plots, as described above for the case in which neither
curve is closed, $d\not=1$, the program then, after scaling each helix in the pair to have approximate length~1 and
scaling the partition discretizing each helix to be a partition of~$[0,1]$, created a common partition of~$[0,1]$ for
the two helices, a nonuniform partition, and discretized each helix by the common partition using cubic splines.
From left to right in the three plots, the common partitions were of size 451, 901, 1051, respectively.
We note as well that
in each case we assumed (correctly) the helices to be defined in the proper directions (see second paragraph of the
Introduction section), thus cutting the times of execution for each case by about half. For each case from left to
right in the three plots, the times of execution were 5.4, 19.4, 39.2 seconds, respectively, with the {\small\bf repeat}
loop in Procedure~$2$ in Section~4 executed 2, 3, 5 times, respectively.
\begin{figure}
\centering
\begin{tabular}{ccc}
\includegraphics[width=0.3\textwidth]{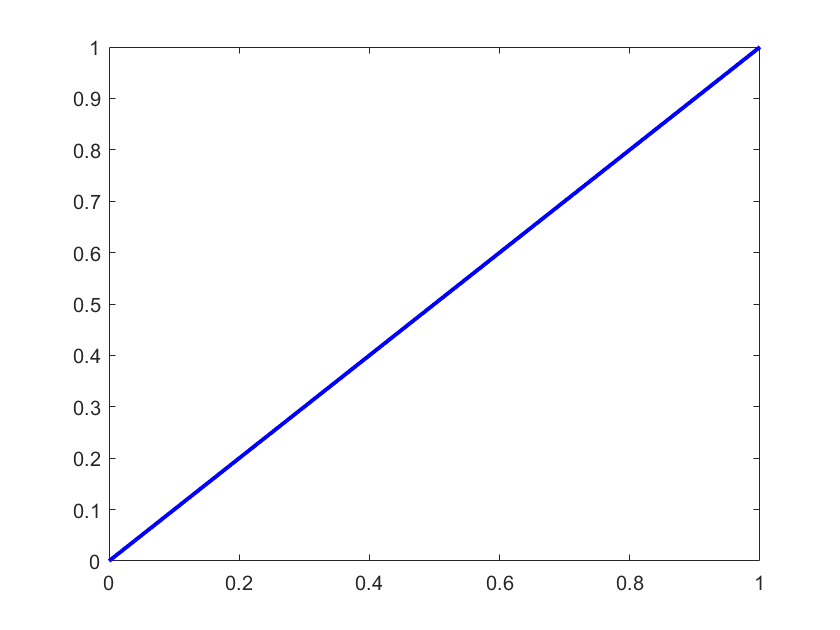}
&
\includegraphics[width=0.3\textwidth]{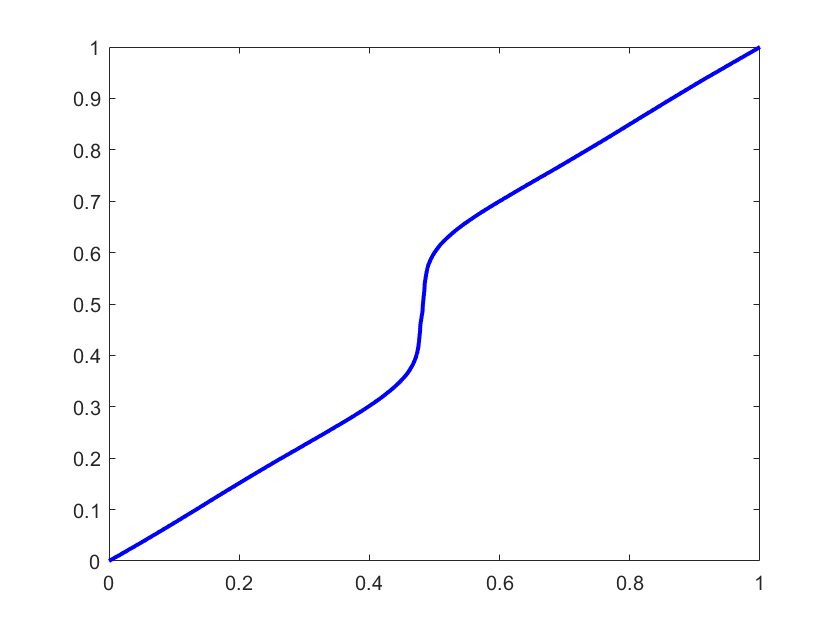}
&
\includegraphics[width=0.3\textwidth]{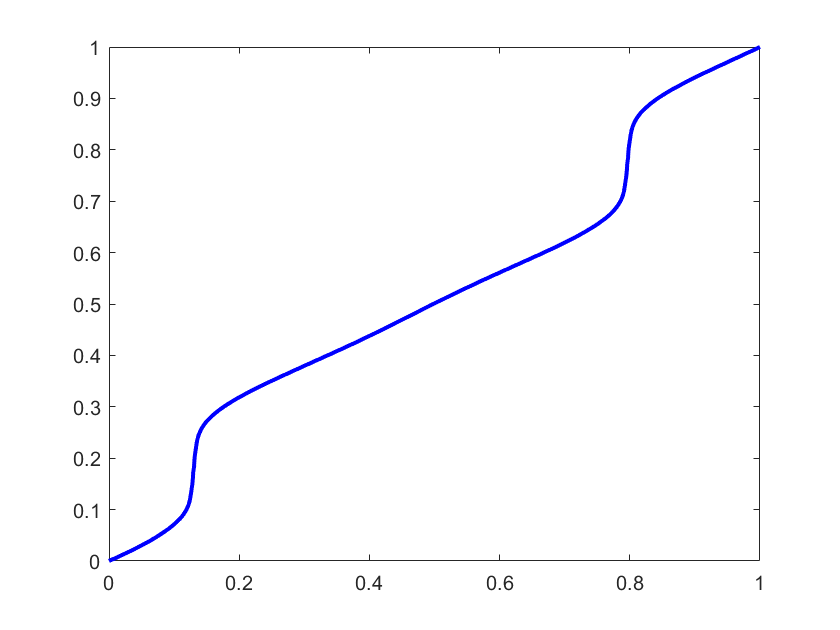}
\end{tabular}
\caption{\label{F:gammas}
Optimal diffeomorphisms for pairs of helices.
}
\end{figure}
\begin{figure}
\centering
\begin{tabular}{ccc}
\includegraphics[width=0.3\textwidth]{curves3.png}
&
\includegraphics[width=0.3\textwidth]{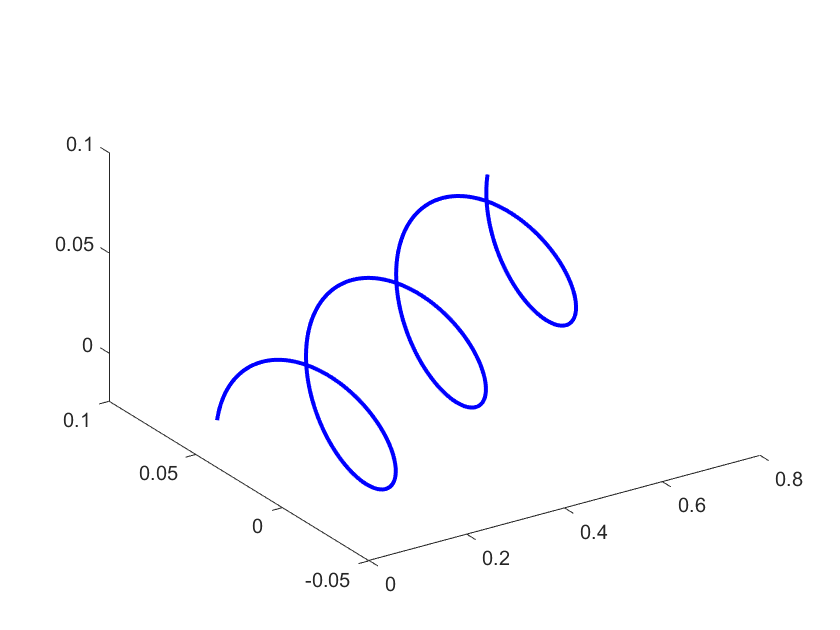}
&
\includegraphics[width=0.3\textwidth]{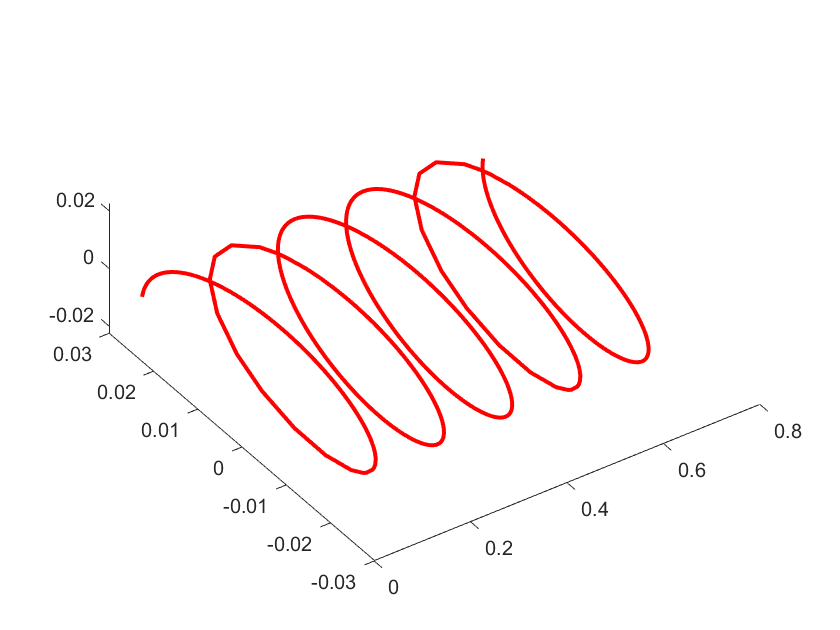}
\end{tabular}
\caption{\label{F:f123best}
Views of first helix of 3 loops and second helix of 5 loops before computation of elastic shape distance and
registration (left), of rotated first helix (middle) and reparametrized second helix (right) after computations.
}

\end{figure}
Plots of the computed optimal diffeomorphisms for each pair of helices from left to right in the three plots in
Figure~\ref{F:curves}, are shown in Figure~\ref{F:gammas}. The computed optimal rotation matrix for the pair in
the leftmost plot in Figure~\ref{F:curves}, was $\left( \begin{smallmatrix} 0 & 0 & 1\\ 1 & 0 & 0\\ 0 & 1 & 0\\
\end{smallmatrix} \right)$. For each of the other two pairs it was almost the same matrix, the entries slightly
different. Finally, Figure~\ref{F:f123best} shows results of the elastic registration of the helix of 3 loops and
the helix of 5 loops. The two helices are shown in the leftmost plot of the figure before any computations took place.
In the middle plot we see the first helix (of 3 loops) after it was rotated with the computed optimal rotation matrix,
its axis of rotation becoming a ray of direction not far from that of the positive $x-$axis. In the rightmost plot
we see the second helix (of 5 loops) after it was reparametrized with the computed optimal diffeomorphism, some of
the consecutive points in its discretization becoming slightly separated, in particular near the end of the first
loop and the beginning of the fifth loop, so that the plot of the helix, which is drawn by joining with line
segments consecutive points in the discretization of the helix, has a slighly flat appearance in these~areas.
\par Finally, we note that we could generate results using spherical ellipsoids similar to the results just
presented for helices. Since such an exercise is tantamount to repeating what has already been done, in its place,
we have opted to use spherical ellipsoids for the purpose of illustrating, if the curves under consideration are
closed, the improvement in execution time that is achieved through the execution of our software package when it
involves the FFT (mostly the execution of Procedure~$3$ in Section~5) as this has the effect of speeding up the
successive computations appearing in Procedure~$2$ in Section~4 of optimal rotation matrices for all starting points
of one of the curves. We also use spherical ellipsoids to illustrate what occurs if the curves under consideration
are not defined in the proper directions.
\begin{figure}
\centering
\begin{tabular}{c}
\includegraphics[width=0.4\textwidth]{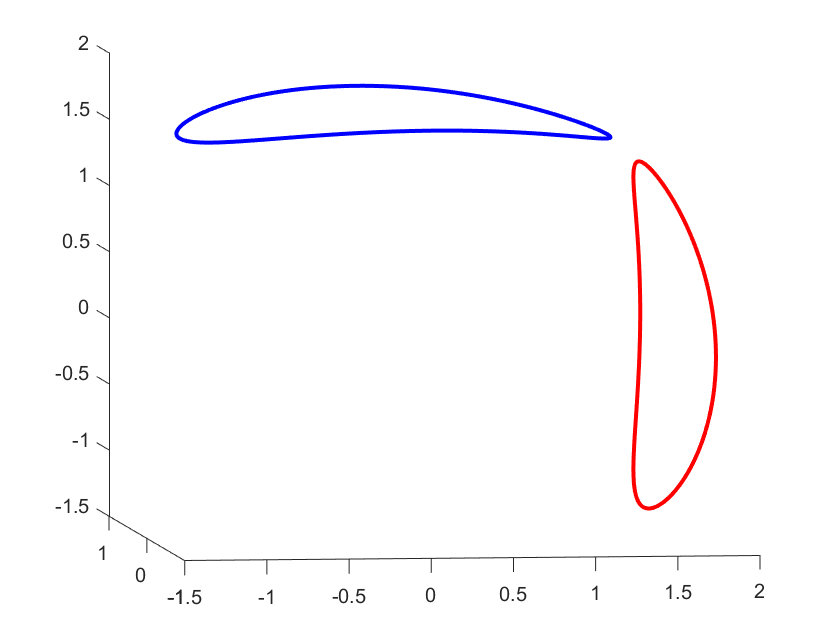}
\end{tabular}
\caption{\label{F:ellips0}
Two spherical ellipsoids, curves in $3-$d space. The positive $z-$axis is the axis of rotation of one
spherical ellipsoid, while the positive $x-$axis is the axis of rotation of the other one. Their shapes
are essentially identical thus the elastic shape distance between them should be essentially~zero.
}
\end{figure}
\par A plot depicting two spherical ellipsoids of essentially the same shape is shown in Figure~\ref{F:ellips0}.
The ellipsoid with the positive $z-$axis as its axis of rotation was considered to be the first curve or ellipsoid
in the plot.  It was obtained by setting $r=2.0$, $a=1.3$, $b=1.0$ in the definition of a spherical ellipsoid above.
The other ellipsoid in the plot has the positive $x-$axis as its axis of rotation and was considered to be the second
curve or ellipsoid in the plot. It was obtained by setting $r=2.0$, $a=1.0$, $b=1.3$ in the definition of a spherical
ellipsoid above, the definition modified in the obvious manner so the the ellipsoid has the positive $x-$axis as
its axis of rotation. The two ellipsoids in the plot were then discretized as described below, and taking into account
that both are closed, the elastic registration of the two ellipsoids and the elastic shape distance
(essentially zero) between them were then successfully computed through the execution of our software package, first without
involving the FFT thus executing mostly Procedure~$2$ in Section~4, and then involving the FFT thus executing
mostly Procedure~$3$ in Section~5. Note that for both procedures the input variable $itop$ was set to~1 as suggested
in Section~4 for spherical ellipsoids.
\smallskip \\
First we discretized the first ellipsoid by a nonuniform partition of $[0,2\pi]$ of size~1001 and the second
ellipsoid by a uniform partition of $[0,2\pi]$ of size~901. As described above for the case in which at least one
curve is closed, the program then, after scaling each ellipsoid to have approximate length~1 and scaling the partition
discretizing each curve to be a partition of~$[0,1]$, created a common partition of~$[0,1]$ for the two ellipsoids,
a uniform partition of size~1001, and discretized each curve by this common partition using cubic splines. The program then
selected the discretization of the first ellipsoid by the uniform partition as the set of starting points of this ellipsoid.
Without involving the FFT, the {\small\bf repeat} loop in Procedure~2 was executed two times, i.e., there was a total of two
iterations for this loop.
The same for the {\small\bf repeat} loop in Procedure~3 when involving the FFT. Without the FFT, the executions of the
KU algorithm for computing successively optimal rotation matrices for all starting points of the first ellipsoid,
took about 0.12 seconds per iteration of the {\small\bf repeat} loop, while the execution of the DP algorithm took
about 6.5 seconds. With the FFT, the execution of the KU2 algorithm, again
for computing successively optimal rotation matrices for all starting points of the first ellipsoid,
took about 0.06 seconds per iteration, while the DP algorithm took about 6.5 seconds.
\smallskip \\
%Replacing above 1001 by 10001 and 901 by 9001, and repeating exactly what was done as described above,
%we then obtained that without the FFT, the executions of the KU algorithm took about 6.00 seconds per iteration
%of the {\small\bf repeat} loop in Procedure~2 (two iterations), while the execution of the DP algorithm took about
%65.10 seconds, and with the FFT, the execution of the KU2 algorithm took about 0.44 seconds per iteration of the
%{\small\bf repeat} loop in Procedure~3 (two iterations), while the DP algorithm took about 65.30 seconds.
%\smallskip \\
Replacing above 1001 by 46001 and 901 by 45001, and repeating exactly what was done as described above,
we then obtained that without the FFT, the executions of the KU algorithm took about 87 seconds per iteration
of the {\small\bf repeat} loop in Procedure~2 (two iterations), while the execution of the DP algorithm took about
291 seconds, and with the FFT, the execution of the KU2 algorithm took about 2 seconds per iteration of the
{\small\bf repeat} loop in Procedure~3 (two iterations), while the DP algorithm took about 291 seconds.
%\smallskip\\
\par
From the two examples above it is clear that computing successively optimal rotation matrices for all starting
points of the first ellipsoid with the FFT is a lot faster than without it. The two examples illustrate as well
the linearity of the DP algorithm and that its execution time appears to be significantly larger than the time
required to compute successively in either Procedure~2 or Procedure~3, optimal rotation matrices for all starting
points of the first ellipsoid. Although the latter may be true when the FFT is used, it is not exactly true otherwise. 
Actually as the size of the discretizations of the curves increases, if the FFT is used, this time becomes
insignificant relative to the execution time of the DP algorithm, but the opposite occurs if it is not.
%although the latter fact seems to obscure the effectiveness of using the FFT as has been described,
\par Finally we reversed the direction of the first ellipsoid in the last example above and as expected obtained
an elastic shape distance between the two ellipsoids different from zero, a distance of 0.195. Using the option
in the program to do the computations in both directions of one of the curves, we then obtained the correct
distance (essentially zero). Of course the execution time of the program doubled.
\\ \smallskip\\ \noindent
%\\ \medskip\\
%\\ \pagebreak\\
{\bf\large Summary}
\\ \smallskip\\
Inspired by Srivastava et al.'s work for computing the elastic registration of two simple curves
in $d-$dimensional Euclidean space, $d$ a positive integer, and thus the associated elastic shape
distance between them, in this paper we have enhanced Srivastava et al.'s work in various ways. First
we have presented a Dynamic Programming algorithm that is linear for computing an optimal diffeomorphism
for the elastic registration of two simple curves in $d-$dimensional space, the computation of the
registration based only on reparametrizations (with diffeomorphisms of the unit interval) of
one of the curves (no rotations), the curves given on input as discrete sets of nodes in the
curves, the numbers of nodes in the curves not necessarily equal, the partitions of the unit
interval discretizing the curves not necessarily uniform. Next we have presented a purely
algebraic justification of the usual algorithm, the Kabsch-Umeyama algorithm, for computing
an optimal rotation matrix for the rigid alignment of two simple curves in $d-$dimensional space,
the curves given on input as discrete sets of nodes in the curves, the same number of nodes in
each curve, the two curves discretized by the same partition of the unit interval, the partition
discretizing the curves not necessarily uniform. Lastly,
with the convention that if one of the curves is closed, the first curve is closed,
%in case one of the curves is closed, assumed to be the first curve,
we have redefined the $L^2$ type distance that is minimized in
Srivastava et al.'s work to allow for the second curve to be reparametrized while the first one
is rotated, the curves again given on input as discrete sets of nodes in the curves, the same
number of nodes in each curve, both curves now discretized by the same partition of the unit
interval (a uniform partition if the first curve is closed). A finite subset of the nodes in the
first curve (possibly all of them, possibly one if neither curve is closed)
is then selected which we interpret to be the set of so-called starting points of the curve,
and the redefined $L^2$ type distance is then minimized with an iterative procedure that
alternates computations of optimal diffeomorphisms (a constant number of them per iteration for
reparametrizing the second curve) with successive computations of optimal rotation matrices
(for rotating the first curve)
for all starting points of the first curve. Carrying out computations this way is not
only more efficient all by itself, but, if both curves are closed, allows applications of the
Fast Fourier Transform (FFT) for computing successively in an even more efficient manner, optimal
rotation matrices for all starting points of the first~curve. We note, results from computations
with the implementation of our methods applied on $3-$dimensional curves of the helix and
spherical ellipsoid kind, have been presented in this paper as~well.

\end{document}